\documentclass[a4paper]{amsart}
\usepackage{txmac}
\usepackage{mypic}
\usepackage{alltt}
\usepackage{myrot}
\usepackage{amssymb}
\usepackage{graphicx}
\usepackage{a4}

\def\mynewtheo#1#2{%
\newtheorem{@#1}[thm]{#2}
\newenvironment{#1}{\begin{@#1}\rm}{\end{@#1}}}

\newtheorem{thm}{Theorem}[section]
\newtheorem{corr}[thm]{Corollary}%[section]
\newtheorem{prp}[thm]{Proposition}

\mynewtheo{exam}{Example}
\mynewtheo{problem}{Problem}
\mynewtheo{question}{Question}
\mynewtheo{rem}{Remark}

\advance\oddsidemargin by -1.7cm
\advance\evensidemargin by -1.7cm
\advance\textwidth by 3.4cm
\advance\textheight by 1.4cm 
\advance\topmargin by -0.7cm 
\begin{document}
\makeatletter

\parskip3mm

\let\sm\setminus
\let\reference\ref
\let\dl\delta
\let\Gm\Gamma
\let\gm\gamma
\let\lm\lambda
\let\sg\sigma
\let\Dl\Delta
\def\cV{{\mathcal V}}
\def\cX{{\mathcal X}}
\def\bQ{{\Bbb Q}}
\def\bN{{\Bbb N}}
\def\bZ{{\Bbb Z}}
\def\sgn{\mathop {\operator@font sgn}\mathord{}}
\def\be{\mathbf{e}}
\let\tl\tilde
\let\es\enspace
\let\ring\r

\long\def\@makecaption#1#2{%
   % \tm
   \vskip 10pt
   {\let\label\@gobble
   \let\ignorespaces\@empty
   \let\myxdef\xdef
   \let\xdef\@empty
   \myxdef\@tempt{#2}%
   \let\xdef\myxdef
   % this work-around is needed for the captions generated by
   % \twofigures (in \twofigures \@tempw, \tempx are set to \@gobble
   % so the so in \xdef\@tempw{..} the \@tempw{..} will disappear
   %\typeout{`#2'}%
   }%
   \ea\@ifempty\ea{\@tempt}{%
   \setbox\@tempboxa\hbox{%
      \fignr#1#2}%
      }{%
   \setbox\@tempboxa\hbox{%
      {\fignr#1:}\capt\ #2}%
      }%
   \ifdim \wd\@tempboxa >\captionwidth {%
      \rightskip=\@captionmargin\leftskip=\@captionmargin
      \unhbox\@tempboxa\par}%
   \else
      \hbox to\captionwidth{\hfil\box\@tempboxa\hfil}%
   \fi}%
\def\fignr{\small\sffamily\bfseries}%
\def\capt{\small\sffamily}%

\let\ig\includegraphics

\newdimen\@captionmargin\@captionmargin2\parindent\relax
\newdimen\captionwidth\captionwidth\hsize\relax

\def\rottab#1#2{%
\expandafter\advance\csname c@table\endcsname by -1\relax
\centerline{%
\rbox{\hsize=\vsize\relax\centerline{\vbox{\setbox1=\hbox{#1}%
\centerline{\mbox{\hbox to \wd1{\hfill\mbox{\vbox{{%
\caption{#2}}}}\hfill}}}%
\vskip9mm
\centerline{
\mbox{\copy1}}}}%
}%
}%
}

\def\br#1{\left\langle\ #1\ \right\rangle}

\def\twofigures#1#2#3#4#5#6{{
% puts two figures beside each other
% #1,#2 - content of figures
% #3,#4 - caption text (WITHOUT \label)
% #5,#6 - the labels (that would appear usually within the \caption)
\let\@tempx\@gobble
\let\@tempw\@gobble
\begin{figure}[h]
\[
\captionwidth0.3\textwidth\relax
\begin{array}{c@{\qquad\quad}c}
#1 & #2 \\
\vbox{\caption{\xdef\@tempw{\@currentlabel}#3}} & 
\vbox{\caption{\xdef\@tempx{\@currentlabel}#4}} 
\end{array}
\]
\end{figure}
\edef\@currentlabel{\@tempw}\label{#5}
\edef\@currentlabel{\@tempx}\label{#6}
}}

\def\mynewtheo#1#2{%
\newtheorem{@#1}{#2}\newenvironment{#1}{\begin{@#1}\rm}{\end{@#1}}}

% \mynewtheo{theorem}{Theorem}
% \mynewtheo{corr}{Corollary}

\def\proof{\@ifnextchar[{\@proof}{\@proof[\unskip]}}
\def\@proof[#1]{\noindent{\it Proof #1.}\enspace}

\title{Mutation and the colored Jones polynomial}
\newbox\@tempboxb
\setbox\@tempboxb=\vbox{\vspace{5mm}

\begin{center}

{\small\it 
BK21 Project, Department of Mathematical Sciences, \\
KAIST, Daejeon 307-701, Korea \\[2mm]
Osaka City University Advanced Mathematical Institute\\
Sugimoto 3-3-138, Sumiyoshi-ku 558-8585, Osaka, Japan\\[2mm]
alexander@stoimenov.net, tanakat@sci.osaka-cu.ac.jp
} \\

\end{center}

\begin{center}
with an appendix by Daniel Matei
\end{center}

% 
% \begin{center}
% DRAFT (\today)
% \end{center}
}
\author{Alexander Stoimenow and Toshifumi Tanaka
\box\@tempboxb
}
%{\bf Toshifumi Tanaka}

\date{Current version: \today\ \ \ First version: Dec 11, 2005}%\date{}

\begin{abstract}
We show examples of knots with the same polynomial invariants and
hyperbolic volumes, with variously coinciding 2-cable polynomials and
colored Jones polynomials, which are not mutants. \\


\end{abstract}

\maketitle

\section{\label{S1}Introduction}

Mutation, introduced by Conway \cite{Conway}, is a procedure 
of turning a knot into another, often different but ``similar'' 
one. This similarity alludes to the circumstance, that most of 
the common (efficiently computable) invariants coincide on mutants, 
and so mutants are difficult to distinguish. A basic exercise in 
skein theory shows that mutants have the same Alexander polynomial 
$\Dl$, and this argument extends to the later obtained Jones $V$, 
HOMFLY (or skein) $P$, BLMH $Q$ and Kauffman $F$ polynomial 
\cite{Jones,HOMFLY,LickMil,BLM,Kauffman}. The cabling formula 
for the Alexander polynomial (see for example \cite[theorem 6.15]
{Lick}) shows also that Alexander polynomials of all satellite knots 
of mutants coincide, and the same was proved by Lickorish and 
Lipson \cite{LicLip} also for the HOMFLY and Kauffman polynomials 
of 2-satellites of mutants. The HOMFLY polynomial applied on a
3-cable can generally distinguish mutants (for example the K-T
and Conway knot; see \S\ref{S32}), but with a calculation effort
that is too large to be considered widely practicable. 

While the Jones polynomial, unlike $\Dl$, was known not to satisfy
a cabling formula (because it distinguishes some cables of knots with
the same polynomial), nontheless Morton and Traczyk \cite{MorTra}
showed that Jones polynomials of all satellites of mutants are
equal. As a follow-up to this result, the question was raised (see
\cite[problem 1.91(2)]{Kirby}; question \reference{q1.1'} below)
whether the converse is true for simple knots. A minorly stronger
paraphrase is:
Is the Jones polynomial of all satellites in fact a \em{universal} 
satellite mutation invariant, i.e. does it distinguish all knots 
which are not mutants or their satellites?

Since the Jones polynomial of all satellites is (equivalent to) 
what is now known as the ``colored Jones polynomial'' (CJP), 
such a universality property relates to some widely studied, 
known or conjectural, features of this invariant. Two important 
recent problems in knot theory, the AJ \cite{G2} and Volume 
conjectures \cite{MM}, assert that the CJP determines the
$A$-polynomial resp. the Gromov norm. Besides, it was proved
(as conjectured by Melvin and Morton) that it determines the
Alexander polynomial \cite{BG,Chmutov,KSA,Vaintrob}, and evidence
is present that it determines the signature function \cite{G}. Latter
problem would be settled, at least for simple knots, by the
aforementioned universality property of the CJP, since it is
known (see \cite{CL}), that mutation preserves all signatures.
Since Ruberman \cite{Ruberman} showed that mutants have equal 
volume in all hyperbolic pieces of the JSJ decomposition, 
universality would imply too (a qualitative version of) the 
Volume conjecture. It is also consistent with the AJ-conjecture 
and a recent result of Tillmann \cite[corollary 3]{Tillmann} 
on coincidence of factors of $A$-polynomials of mutants. Note
also that the Volume conjecture\footnote{The second author
\cite{T} has shown that it is in fact sufficient the Volume 
conjecture to hold for doubled knots.}, as well as the AJ-conjecture
\cite{DG}, in turn imply that the colored Jones polynomial, 
and hence Vassiliev invariants \cite{BN}, detect the unknot.

The main motivation for this paper is to answer Question \ref{q1.1'}.

\begin{thm}\label{rg}
There exist infinitely many pairs of (simple) hyperbolic knots with
equal CJP, which are not mutants.
\end{thm}

Below we will show constructions of knots with equal colored Jones 
polynomial. The exclusion of mutation in theorem \ref{rg} bases on a 
study Vassiliev invariants obtained from (cabled) knot polynomials.
Apart from providing such examples, we will examine closer various 
other criteria for mutation. We will use also a different concept,
representations of the fundamental group of the double branched 
cover (contributed by Daniel Matei in an appendix to our paper). 
Either of the Vassiliev invariant and the group representations 
approach may be more useful than the other, as we show by examples. 
In contrast, we give in \S\ref{S4} also instances where CJP excludes 
mutation, but other invariants (knot polynomials, hyperbolic volume) 
fail. Some of our arguments are followed by several remarks. These 
try to address the relation of our examples to the AJ- and Volume
conjecture, as well as combinations of mutation criteria for which
we do not know if (non-)distinction phenomena occur. The conclusions 
our work allows us to make can be summed up like this (see remark 
\ref{r35}, examples \ref{x41} and \ref{x56}, and theorem \ref{thach}):

{\parskip\z@
\begin{thm}\label{th1.2}
\def\myitem#1{\noindent\hangafter1\hangindent0.7cm\hbox to
0.7cm{\rm #1\hfil}\kern0pt}
\mbox{}

\myitem{(1)}\ry{1.4em}%
The CJP does not determine the HOMFLY, Kauffman polynomial or
their 2-satellites, or the fundamental group of the double
branched cover.

\myitem{(2)}
The CJP of hyperbolic knots is not determined by hyperbolic volume,
the double cover, the HOMFLY, Kauffman polynomial and 2-cable HOMFLY
polynomial, even when all of them are taken in combination.

\myitem{(3)}
The property two hyperbolic knots to be mutants is not determined
by the coincidence (even in combination) of hyperbolic volume, CJP,
HOMFLY, Kauffman polynomial and either (a) their 2-satellites or
(b) the double cover.
\end{thm}
}

A brief outline of the paper is as follows. After \S\ref{S2},
containing some preliminaries on the Colored Jones polynomial,
we will start in \S\ref{S3} with some examples to prove theorem 
\ref{rg}. These, and many of the following, examples grew out of 
the first author's attempt to determine the mutations among low 
crossing knots \cite{mut}. In  \S\ref{S4} we study some pairs 
consisting of a knot and its mirror image, and then further refine 
the construction for and proof of Theorem \ref{rg} to adapt it to 
such pairs (see Theorem \ref{thach}). We have, however, also cases, 
where the polynomial invariants fail. These examples are shown in 
\S\ref{S5}. For such knots the exclusion of mutation was extremely 
difficult, and we were assisted by Daniel Matei. He explains his 
calculation in the appendix of the paper.

%########################### attached ###################################

\section{The colored Jones polynomials of knots\label{S2}}

\noindent First we set up and clarify some terminology concerning
satellites.

Let $L$ be a link embedded in the solid torus $T=S^1\times D^2$. If
we embed $T$ in $S^3$ so that its core $S^1\times \{0\}$ represents
a knot $K$, then we call the resulting embedding $K'$ of $L$ the
\em{satellite of (companion) $K$ with pattern $L$}.
The satellite is defined up to
the choice of framing of $T$. The \em{algebraic intersection number}
$d_a(L)$ of $L$ is its absolute homology class in $H_1(T)=\bZ$.
The \em{geometric intersection number} $d_g(L)$ of $L$ is the
smallest number of transverse intersection points of $L$ with
a meridional disk $D$ of $T$. Clearly $d_g(L)\ge d_a(L)$ for any
$L$. If $K'$ is a satellite with pattern $L$, we call the number
$d_g(L)$ also the \em{degree} of $K'$. If $d_a(L)=d_g(L)$ (i.e.\ all
the intersections of $D$ with $L$ can be made so that $L$ points
in the same direction w.r.t.\ $D$), then we call the satellite $K'$
also a \em{cable} of $K$. A cable $K'$ of $K$ is called \em{connected}
if $K'$ (or the pattern $L$) is a knot, and \em{standard}, if
$L=S^1\times \{x_1,\dots,x_n\}$, with $x_i$ being distinct
points in $D^2$. A satellite/cable of degree $n$ is simply called
an $n$-satellite or $n$-cable. If $L$ is one of the components
of the Whitehead link, and $T$ the complement of the other component,
then we call $K'$ a \em{Whitehead double} of $K$.

\noindent
Let $\mathbb{Z}[A,A^{-1}]$ be the Laurent polynomial ring in one
indeterminate $A$ with coefficients in the ring of integers. We put
$q=a^2=A^4$ (so that $q=1/t$ for the variable $t$ of the Jones
polynomial standardly used), and set $\{n\}=a^{n}-a^{-n}$, $[n]=
\{n\}/\{1\}$ and $[n]!=[1][2]\ldots[n]$.

\noindent
The {\it Kauffman bracket skein module} $\mathcal{K}(M)$ of an oriented 
3-manifold $M$ is the quotient of the free $\mathbb{Z}[A,A^{-1}]$-module
generated by the set of ambient isotopy classes of framed links in $M$, 
by the following {\it Kauffman relations}:
\vspace{1mm}
% \begin{enumerate}
\[
<L\coprod
\raisebox{-0.15cm}{\ig[width=0.5cm,height=0.5cm]{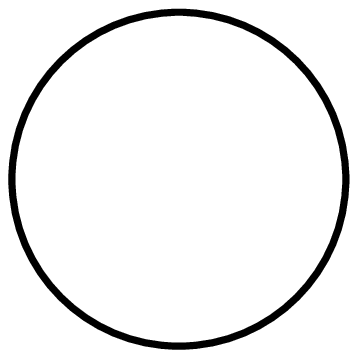}}>\,
=-[2]<L>,
\qquad <
\raisebox{-0.2cm}{\ig[width=0.5cm,height=0.5cm]{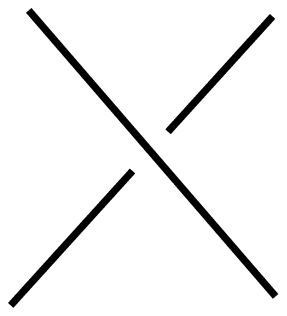}}
>\,=A<
\raisebox{-0.2cm}{\ig[width=0.5cm,height=0.5cm]{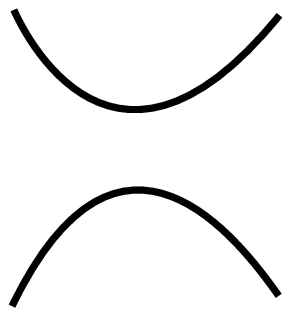}}
>+A^{-1}<
\raisebox{-0.2cm}{\ig[width=0.5cm,height=0.5cm]{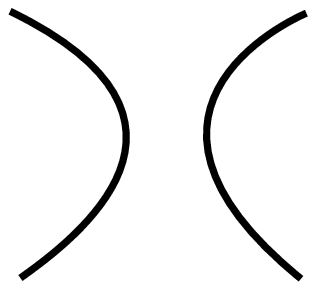}}
>\,.
\]
% \end{enumerate}

\noindent
We know that $\mathcal{K}(S^3)$ is identified with $\mathbb{Z}[A,A^{-1}]$.
The {\it Kauffman bracket} $<L>\in\mathbb{Z}[A,A^{-1}]$ of a framed link $L$ 
in $S^3$ is defined by the image of the isomorphism $\mathcal{K}(S^3)\to
\mathbb{Z}[A,A^{-1}]$ which takes the empty link $\varnothing$ to $1$.\\

\noindent
The skein module of the solid torus $T=S^1\times D^2$ is $\mathbb{Z}[A,A^{-1}][z]$. 
Here $z$ is given by the framed link $S^1 \times I$, where $I$ is a small arc in the 
interior of $D^2$, and $z^n$ means $n$ parallel copies of $z$. Put $\mathcal{B}=
\mathcal{K}(T)$. There is a basis $\{\be_{i}\}_{i\ge 0}$ for 
$\mathcal{B}$ which is defined recursively by \\[-0.6mm]
\begin{equation}\label{En}
\be_{0}=1, \be_{1}=z, \be_{i}=z\be_{i-1}- \be_{i-2}\,.
\end{equation}
\\[-2mm]
\noindent
Let $K$ be a knot in $S^3$. We assume that $K$ is equipped with the zero framing. 
We define a $\mathbb{Z}[A,A^{-1}]$-linear map $<\ >_{K}:\mathcal{B}\rightarrow 
\mathbb{Z}[A,A^{-1}]$ for $K$ by cabling $K$ and taking the Kauffman bracket.
The {\it $N$-colored Jones polynomial} of a knot $K$ is defined as the Kauffman 
bracket of $K$ cabled by $(-1)^{N-1} \be_{N-1}$:\\[-0.6mm]
\[
J'_{K}(N;A)=(-1)^{N-1}<\be_{N-1}>_{K}\,.
\]
\\[-2mm]
It is then normalized
\[
J_{K}(N;A)=\frac{J'_{K}(N;A)}{J'_{\bigcirc}(N;A)}\,,
\]
\\[-2mm]
so that it takes value 1 on the unknot.

\noindent
This invariant is the quantum invariant corresponding to 
the $N+1$-dimensional representation of $sl_2$. As in \S\ref{S1},
we continue using the abbreviation CJP.

\noindent
Now we explain the graphical calculus of G. Masbaum
and P.\ Vogel \cite{MV} (or see also \cite{BHMV}).\\

\noindent
Let $\mathbb{Q}(A)$ be the field generated by the indeterminate $A$ over 
the rational numbers $\mathbb{Q}$. Framed $(n,n)$-tangles with Kauffman relations 
generate a finite-dimensional associative algebra $T_{n}$ over $\mathbb{Q}(A)$, 
which is called the {\it Temperley-Lieb algebra on n-strings}.
$T_{n}$ is generated by the following elements.

% \begin{figure}[h]
\[
\ig[width=9cm,height=2.8cm]{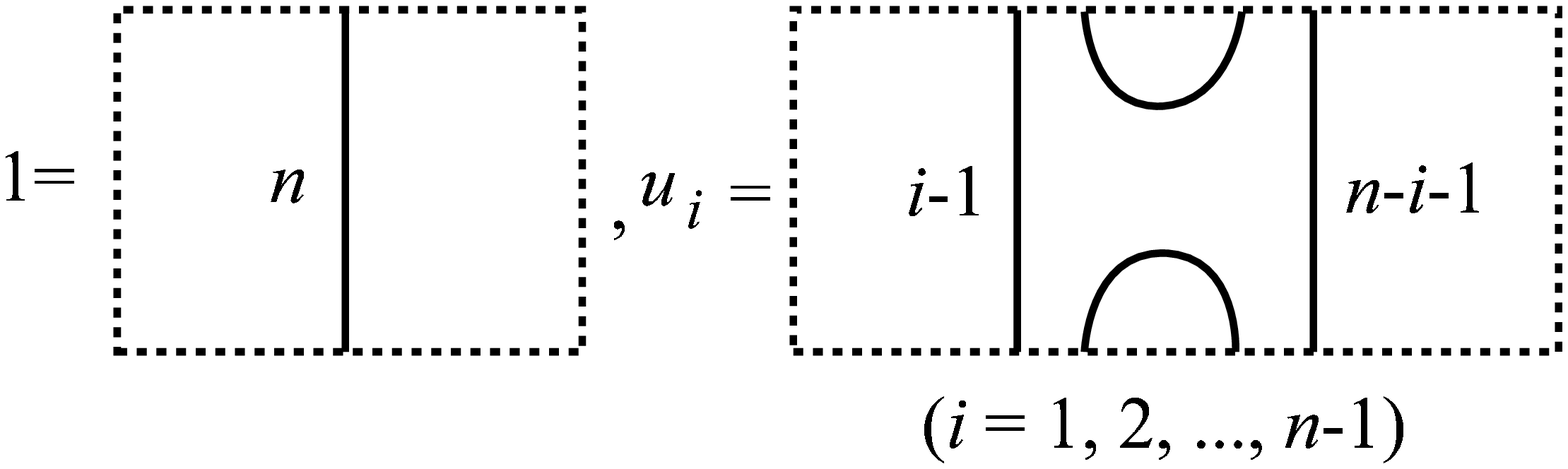}
\]
%		\caption{}
%	\label{}
% \end{figure}
% need eps2eps -dEPSCrop or -dEPSFitPage for Pic6.eps compression
% because BoundingBox had negative value

\noindent
An integer beside an arc signifies $n$ copies of the arc all parallel in the plane. 
There is a {\it trace map} $g:T_{n}\rightarrow \mathcal{B}$ given by mapping a tangle 
with square to the diagram in the annulus obtained by identifying the upper and lower 
edges of the diagram. If we put $d_{n}=g(f^{(n)})$, where $f^{(n)}$ is the {\it 
Jones-Wenzl idempotent} in $T_{n}$, then it is well known that $\{d_{j}\}_{j\ge 0}$
satisfies the recurrence formula for $\{\be_{i}\}_{i\ge 0}$. So the basis 
$\{\be_{i}\}_{i\ge 0}$ can be defined by using the Jones-Wenzl idempotent 
and the trace map. \\

\noindent
An {\it admissibly colored framed (or ribbon) trivalent graph}
is defined as follows.
A $color$ is just a nonnegative integer. A triple of colors $(a,b,c)$
is {\it admissible} if it satisfies the following conditions: \\[-5mm]
\begin{itemize}
\item $a+b+c$ is even, and \\[-2mm]
\item $a+b \ge c \ge |a-b|$, as well as the analogous two other
inequalities obtained by permuting $a,b,c$.
\end{itemize}
Let $D$ be a planar diagram of a ribbon trivalent graph. An
admissible coloring of $D$ is an assignment of colors to the edges
of $D$ so that at each vertex, the three colors meeting there form
an admissible triple.\\

\begin{figure}[h]
\ig[width=5.6cm,height=4cm]{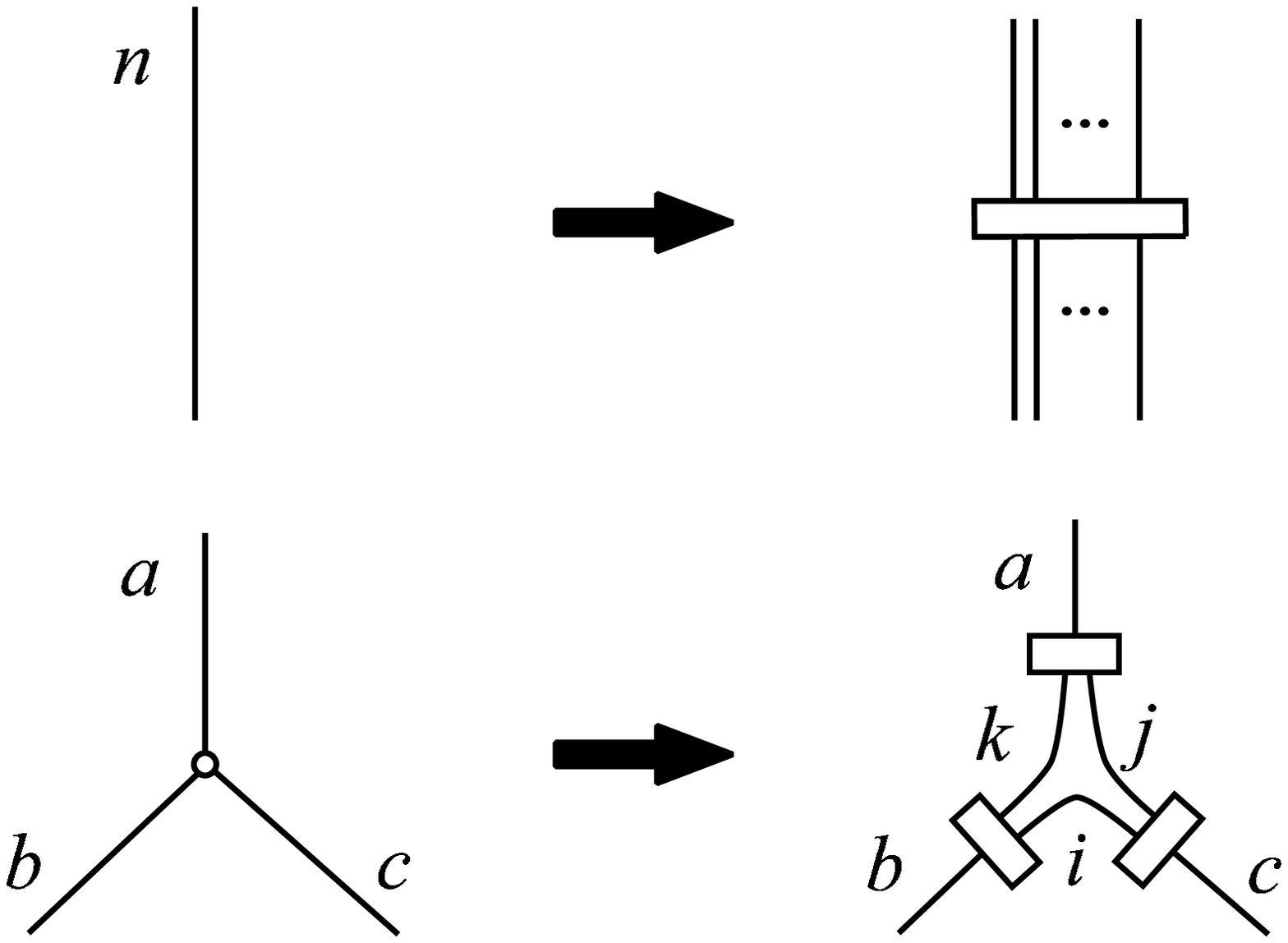}
		\caption{}
	\label{Pic7}
\end{figure}

\noindent
The Kauffman bracket of an admissibly colored framed trivalent graph $D$ is defined 
to be the Kauffman bracket of the expansion of $D$ obtained as follows. The expansion 
of an edge colored by $n$ consists of $n$ parallel strands with a copy of the {\it 
Jones-Wenzl idempotent} $f^{(n)}$ inserted and the colored vertices are expanded as 
in Figure \ref{Pic7}. 

\noindent
The Jones-Wenzl idempotent $f^{(n)}$ is characterized by the property
that $u_{i}f^{(n)}= f^{(n)}u_{i}=0$ for each $u_{i}$ $(1\le i\le n-1)$
and $(f^{(n)})^{2}=f^{(n)}$. The idempotent is represented by a
little box. The triple $(a,b,c)$ is admissible. The {\it 
internal colors} $i,j,k$ are defined by 
\begin{equation}\label{sv}
i=(b+c-a)/2,\es\, j=(a+c-b)/2,\es \mbox{and}\es k=(a+b-c)/2\,.
\end{equation} 

\noindent
Now we can regard $<\be_{n}>_{K}$ as a planar diagram of a
ribbon trivalent graph
by writing $n$ beneath $K$ and inserting a little box into $K$. \\

To simplify language, let us call a \em{satellite of the Jones
polynomial} an invariant of knots $K$ obtained by evaluating the
(ordinary) Jones polynomial on a satellite of $K$. (In other words,
we dualize the satellite operation on the level of invariants.)

This way, the colored Jones polynomial $J_K$ of a knot $K$ can be
considered as a sequence of polynomials $\{J_K(N;A)\}_{N\in \bN}$
in $A$, associated to natural numbers $N$. These polynomials are
obtained by evaluating the Jones polynomial on a cable of $K$,
decorated with the Jones-Wenzl idempotent. So $J_{\star}(N;A)$
are linear combinations of satellites of the Jones polynomial. 
For $N=1$ we have the Jones polynomial $V$ itself.

The Kauffman bracket skein module of a decoration is obtained
from crossing-free connections of the $N$ strands on top and
bottom, and decoration with all crossing-free connections except
the parallel one ($N$ parallel strands) gives the Jones polynomial
of a lower order cable. Now the elements $\be_N$ of \eqref{En}
have a non-zero coefficient in the parallel crossing-free connection
$z^N$. This observation establishes that the information contained
in all satellites of the Jones polynomial is equivalent to the
colored Jones polynomial. It also shows that, for the lowest $N$
distinguishing two knots, it is equivalent whether we talk of the
$N$-cabled or $N$-colored polynomial.

\noindent
We put, with $i,j,k$ as in \eqref{sv},
\[
<k>\,=(-1)^{k}[k+1]\,,\qquad <a,b,c>\,=
(-1)^{i+j+k}\frac{[i+j+k+1]![i]![j]![k]!}{[i+j]![j+k]![i+k]!}\,.
\]

\noindent
Then G. Masbaum and P. Vogel showed the following ``fusion''
formula \cite{MV}:
\vspace{2mm}
% 
% \tm
\begin{equation}\label{MV}
% \displaystyle
\raisebox{-0.9cm}{\mbox{\ig[width=2cm,height=2cm]{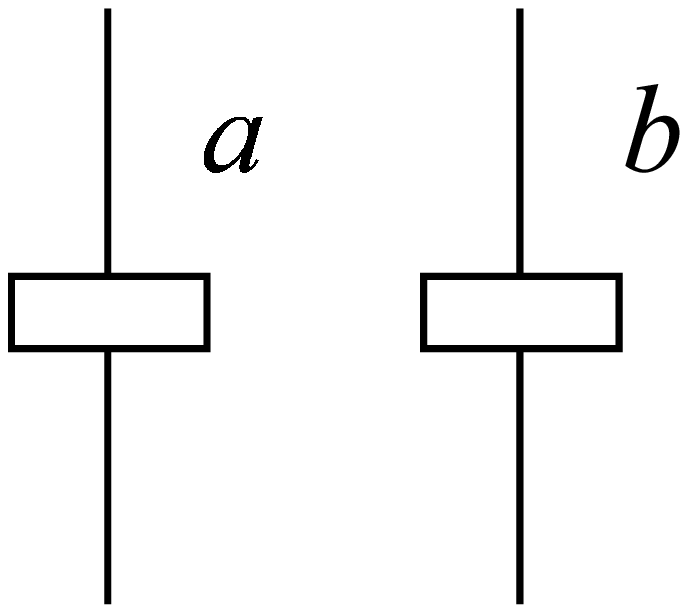}}}
\quad
=\quad\sum_{c}\frac{<c>}{<a,b,c>}\quad \raisebox{-0.9cm}{\mbox{
  \ig[width=0.8cm,height=2cm]{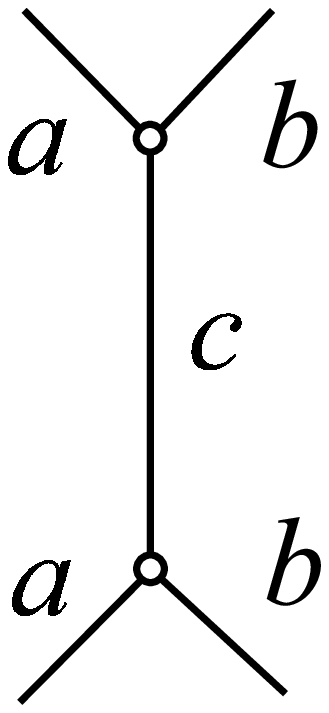}}}\,,
\end{equation}\\[2mm]
where the sum runs over $c$ so that $(a,b,c)$ is admissible.

This formula will enable us to show that certain CJP
are equal, though it requires to find a convenient
presentation of the knots (and the concrete evaluation
of the polynomials remains highly uninviting).

% \section{Examples of non-mutants\label{S3}}

\section{Mutations and invariance of the
Colored Jones polynomial\label{S3}}

\subsection{Initial examples}

\noindent
A \em{mutation} (in Conway's original\footnote{There are now 
various extensions of this concept, also to 3-manifolds; see e.g.
\cite{CL,Rong}.} sense \cite{Conway}) is the following operation.
Consider a knot being formed from two tangles $T_{1}$ and $T_{2}$.
(A tangle is understood here to consist of two strings.) Cut the 
knot open along the endpoints on each of the four strings coming 
out of $T_{2}$. Then rotate $T_{2}$ by $\pi$ along some of the 3
axes - horizontal in, vertical in, or perpendicular to the projection
plane. This maps the tangle ends onto each other. Finally, glue
the strings of $T_{1,2}$ back together (possibly altering orientation
of all strings in $T_{2}$). Two knots $K_1$ and $K_2$ are \em{mutants}
if they can be obtained from each other by a sequence of mutations.
\\

A knot which is the connected sum of two non-trivial knots is called
\em{composite}. We say that a knot $K'$ is a \em{satellite knot} if 
it is a satellite of degree at least 2 around some knot $K$. (Note
that a knot $K'$, which is a satellite of degree $1$ around $K$, is
either equal to $K$, or the connected sum of $K$ with some other knot.)
A knot $K$ is called \em{simple} if it is not a composite or satellite
knot, or in other terms, if its complement has no incompressible tori.
By Thurston's work \cite{Thurston1,Thurston2}, such a knot is
hyperbolic or a torus knot.

\noindent
The central aim of this paper is to answer the following question. \\

% \noindent
\begin{question}\label{q1.1'}
({\rm Problem} 1.91(2) in \cite{Kirby}). Let $K$ be a simple
unoriented knot. Are there any knots other than mutants of $K$,
which cannot be distinguished from $K$ by the Jones polynomial
and all its satellites?
\end{question}

\begin{figure}[h] 
	\ig[width=8.0cm,height=3cm]{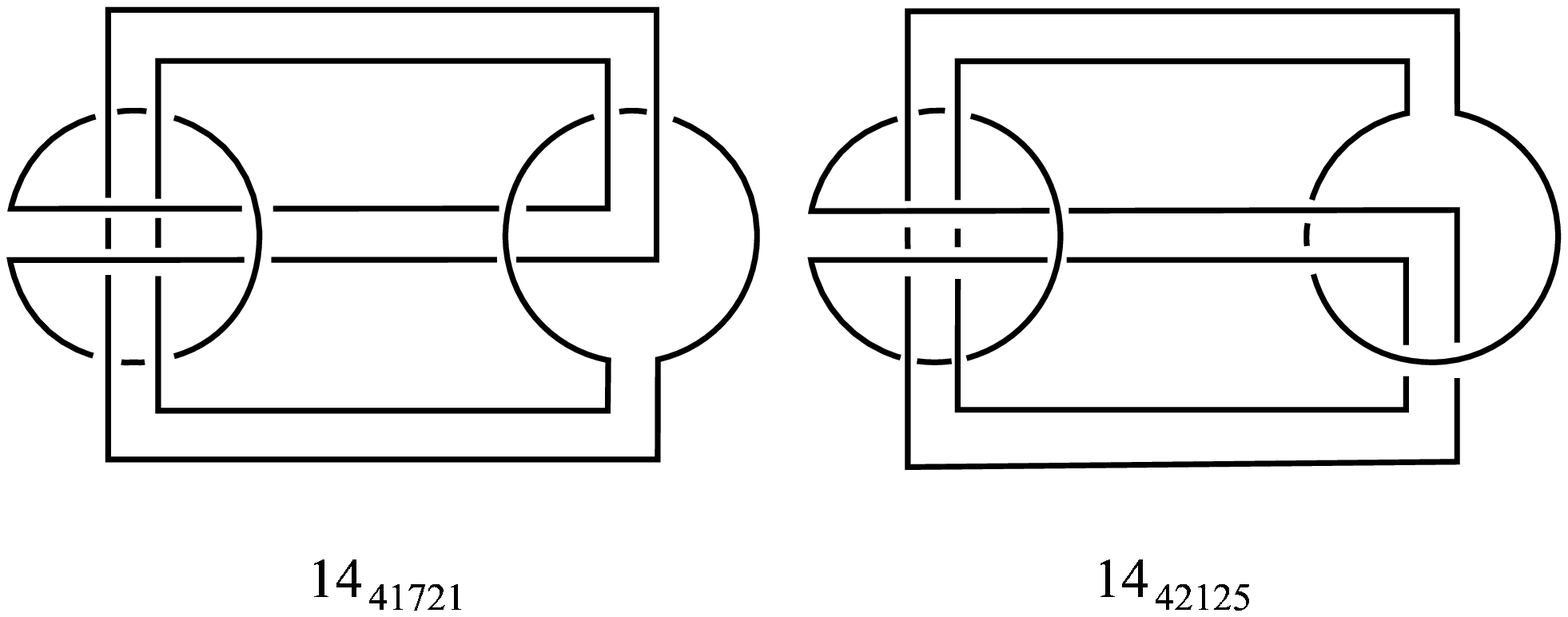}
	\caption{}
	\label{Pic1}
\end{figure}

\noindent
The knots $14_{41721}$ and $14_{42125}$ from \cite{HT}
are depicted in Figure \ref{Pic1}. (In the following we
use the knot ordering in \cite{HT} throughout, but we append
non-alternating knots of given crossing number after the
alternating ones, and we do not take care of mirror images,
unless it is relevant.)
% 
% 
% \noindent
These are ribbon knots, and have trivial Alexander polynomial.
$14_{41721}$ and $14_{42125}$ have the same HOMFLY and Kauffman
polynomial invariants, the same 2-cable HOMFLY polynomial, and the
same hyperbolic volume. For this pair of knots, we show the 
following.

\begin{prp}\label{th11}
$14_{41721}$ and $14_{42125}$ are not mutants.
\end{prp}

\begin{prp}\label{th12}
$14_{41721}$ and $14_{42125}$ have the same colored Jones polynomial.
\end{prp}

Propositions \ref{th11} and \ref{th12} combinedly give the
first examples answering Question \ref{q1.1'}.

A remark on methodology is in place.
Note that, by skein theory, for given strand orientation we
need to calculate the HOMFLY polynomials $P$ of two degree-2
satellites to ensure with certainty that the polynomials of all
degree-2 satellites (of the given strand orientation) coincide.
(In opposition, for the Jones polynomial one satellite for an
arbitrarily chosen orientation suffices in any degree.) For 2-cable
we used the knots obtained by blackboard framing from the diagram $D$
in the table of \cite{HT} and its mirror image $!D$, both 2-cables
having one negative half-twist crossing. For Whitehead double we
used the (framed) satellites of $D$ and $!D$ with a positive clasp.
Similar is the case of the Kauffman $F$ polynomial, though there
it is not necessary to distinguish strand orientation, and so we
consider only the 2-cables.

\begin{table}[h]
{\small%\scriptsize%tiny%footnotesize
\begin{verbatim}

14.   41721:
0 22
 -10   8                  9     57    142    174     98      3    -26    -10      1      1
 -12   8          18    -40   -496  -1284  -1588   -984   -122    246    132    -10    -16
 -12   8        -138    105   2229   5257   5895   3693    613  -1207   -843      0     92
 -12   8         449   -253  -6064 -12412 -11763  -7097  -1346   3258   2546     88   -238
 -12   8        -744    449  10297  18323  13979   7797   1633  -5010  -4098   -182    310
 -12   8         680   -470 -11184 -17574 -10362  -5112  -1181   4587   3846    156   -212
 -12   8        -354    277   7919  11167   4871   2016    516  -2576  -2210    -65     77
 -12   8         104    -90  -3680  -4724  -1440   -466   -132    892    787     13    -14
 -12   8         -16     15   1109   1313    257     58     18   -185   -169     -1      1
 -12   4           1     -1   -208   -230    -25     -3     -1     21     20
  -8   4                        22     23      1      0      0     -1     -1
  -8  -6                        -1     -1

14.   42125:
0 22
 -10   8                  9     57    142    174     98      3    -26    -10      1      1
 -12   8          16    -56   -550  -1380  -1672   -984    -38    342    186      6    -14
 -12   8        -136    149   2451   5745   6371   3693    137  -1695  -1065    -44     90
 -12   8         449   -295  -6414 -13392 -12841  -7097   -268   4238   2896    130   -238
 -12   8        -744    465  10571  19331  15237   7797    375  -6018  -4372   -198    310
 -12   8         680   -472 -11298 -18148 -11176  -5112   -367   5161   3960    158   -212
 -12   8        -354    277   7943  11349   5163   2016    224  -2758  -2234    -65     77
 -12   8         104    -90  -3682  -4754  -1494   -466    -78    922    789     13    -14
 -12   8         -16     15   1109   1315    261     58     14   -187   -169     -1      1
 -12   4           1     -1   -208   -230    -25     -3     -1     21     20
  -8   4                        22     23      1      0      0     -1     -1
  -8  -6                        -1     -1

\end{verbatim}
}
\caption{The polynomials $P_{w,+}$ for $14_{41721}$ and $14_{42125}$.}
\label{figP}
\end{table}

\noindent
{\it Proof of Proposition \reference{th11}.} 
Let $P_{w,+}$ be the evaluation of the HOMFLY $P$ polynomial
on the 0-framed Whitehead double with positive clasp.
The non-mutant status of the knots could be shown by calculation
of the $P_{w,+}$ polynomials, shown in Table \reference{figP}.
The convention
for the HOMFLY polynomial is so that the skein relation has,
as in \cite{LickMil}, the form
\[
l\,P\bigl(
\diag{5mm}{1}{1}{
\picmultivecline{0.18 1 -1.0 0}{1 0}{0 1}
\picmultivecline{0.18 1 -1.0 0}{0 0}{1 1}
}
\bigr)\,+\,
l^{-1} \,P\bigl(
\diag{5mm}{1}{1}{
\picmultivecline{0.18 1 -1 0}{0 0}{1 1}
\picmultivecline{0.18 1 -1 0}{1 0}{0 1}
}
\bigr)\,=\,
-m\,P\bigl(
\diag{5mm}{1}{1}{
\piccirclevecarc{1.35 0.5}{0.7}{-230 -130}
\piccirclevecarc{-0.35 0.5}{0.7}{310 50}
}
\bigr)\,,
\]
and the unknot has unit polynomial. %It is as in \cite{LickMil}. , but
% with $l$ and $l^{-1}$ interchanged (so the right-hand trefoil has
% a polynomial with positive powers of $l$).
The line below the knot
denotes the minimal and maximal degree in the $m$-variable. The
coefficients in (increasing even powers of) $m$, which are polynomials
in $l$, follow line by line, with the minimal and maximal degree in $l$
recorded first and followed by the coefficients of (even powers of)
$l$. \qed

\begin{rem}\label{rK}
Calculation showed that
the 2-cable Kauffman polynomial also distinguishes the two knots.
By an easy skein argument, so will then the Whitehead double
Kauffman polynomials (with some framing at least).
\end{rem}

\noindent
{\it Proof of Proposition \reference{th12}.}
By making repeatedly use of the formula \eqref{MV}, we have 

\noindent
$J_{14_{41721}}(A; N)$
$=\displaystyle \sum_{k_{1}=0}^{N}\frac{<2k_{1}>}{<N,N,2k_{1}>}
\br{\raisebox{-1cm}{\ig[width=4cm,height=2.5cm]{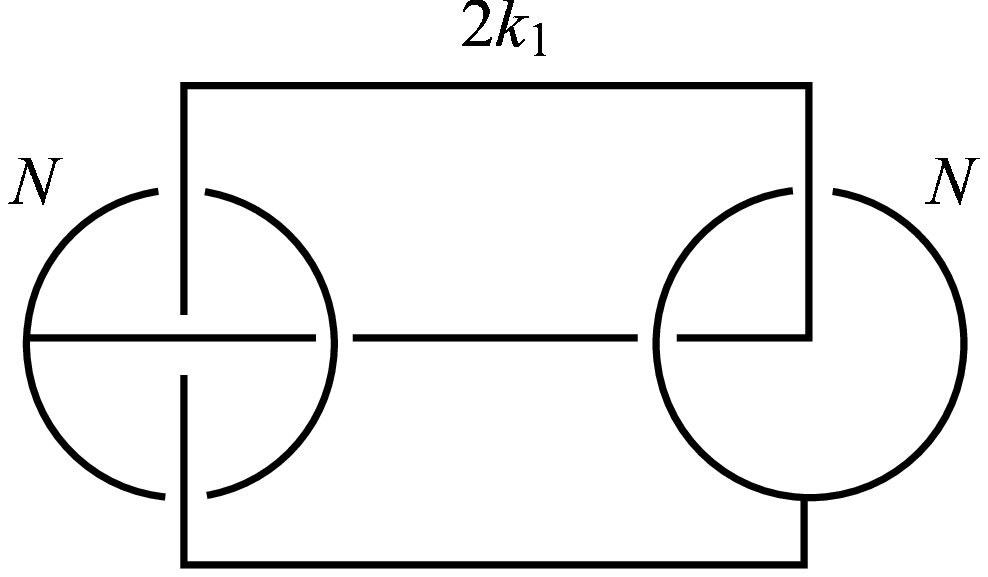}}}$\\

\noindent
$=\displaystyle \sum_{k_{1}=0}^{N}\sum_{k_{2}=0}^{2k_{1}}\frac{<2k_{1}>}{<N,N,2k_{1}>}
\frac{<2k_{2}>}{<2k_{1},2k_{1},2k_{2}>}
\br{\raisebox{-1cm}{\ig[width=4cm,height=2.5cm]{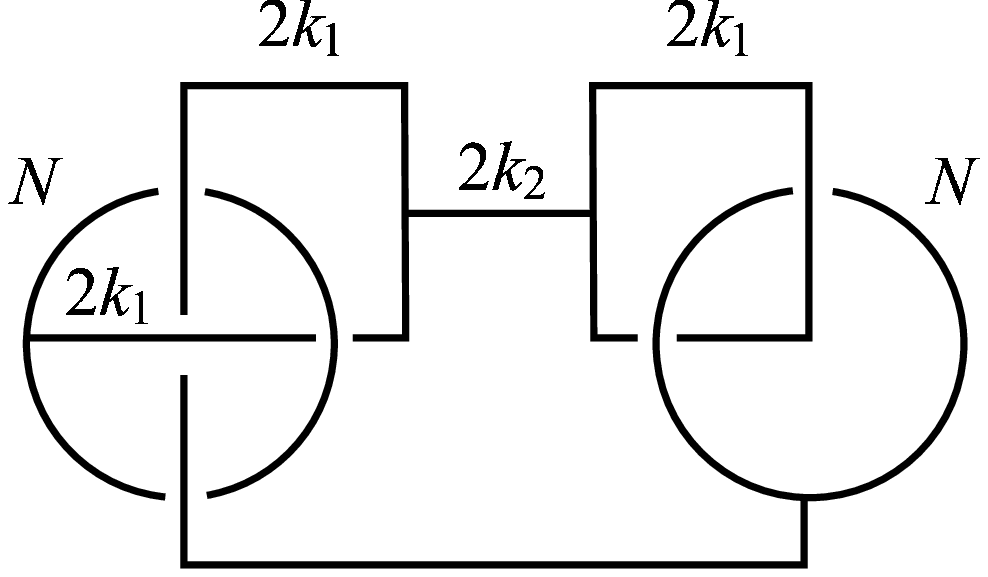}}}$\\

\noindent
$=\displaystyle \sum_{k_{1}=0}^{N}\sum_{k_{2}=0}^{2k_{1}}\sum_{k_{3}=0}^{k_{1}+k_{2}}
\frac{<2k_{1}>}{<N,N,2k_{1}>}\frac{<2k_{2}>}{<2k_{1},2k_{1},2k_{2}>}
\frac{<2k_{3}>}{<2k_{1},2k_{2},2k_{3}>}\ %\times\\ \times
\br{\raisebox{-1cm}{\ig[width=4cm,height=2.5cm]{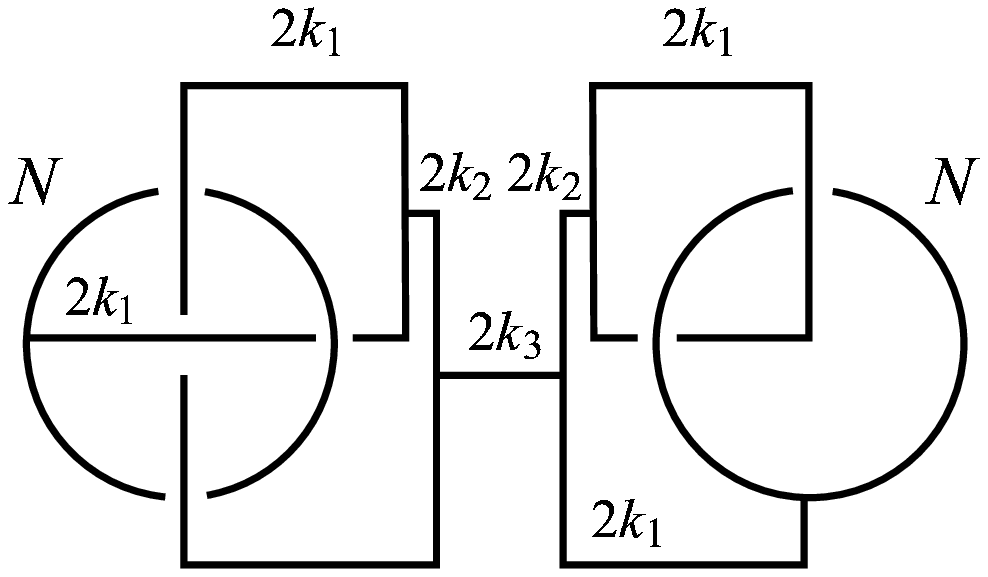}}}$, \\

\noindent
where the third sum runs over all $k_3$ such that $(2k_{1}, 2k_{2},
2k_{3})$ is admissible. We note that $k_{3}$ equals zero by a property
of the Jones-Wenzl idempotent. This implies that the term in the sum
vanishes if $k_{1}\neq k_{2}$. So we have

\noindent
$J_{14_{41721}}(A; N)$
\noindent
$=\displaystyle \sum_{k_{1}=0}^{N}\frac{<2k_{1}>}{<N,N,2k_{1}>}
\frac{<2k_{1}>}{<2k_{1},2k_{1},2k_{1}>}
\br{\raisebox{-1cm}{
\ig[width=2.0cm,height=2.5cm]{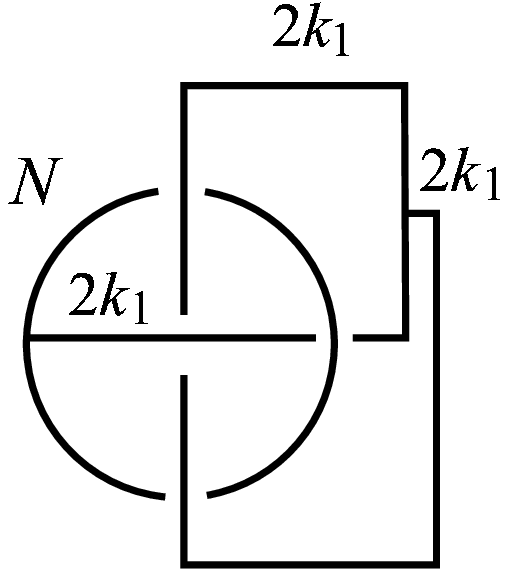}
}}\ 
% \times\\ \times
\br{
\raisebox{-1cm}{\ig[width=1.8cm,height=2.5cm]{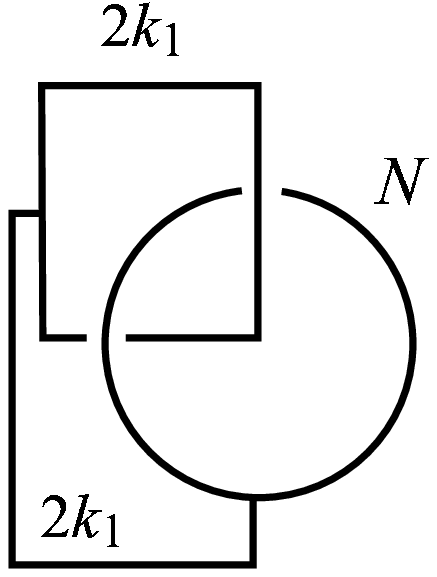}}
}$. \\

\noindent
By using the same method, we have

\noindent
$J_{14_{42125}}(A; N)$
\noindent
$=\displaystyle \sum_{k_{1}=0}^{N}\frac{<2k_{1}>}{<N,N,2k_{1}>}
\frac{<2k_{1}>}{<2k_{1},2k_{1},2k_{1}>}
\br{
\raisebox{-1cm}{\ig[width=2.0cm,height=2.5cm]{eps/Pic13.eps}
}}
\ %\times\\ \times
\br{
\raisebox{-1cm}{\ig[width=1.8cm,height=2.5cm]{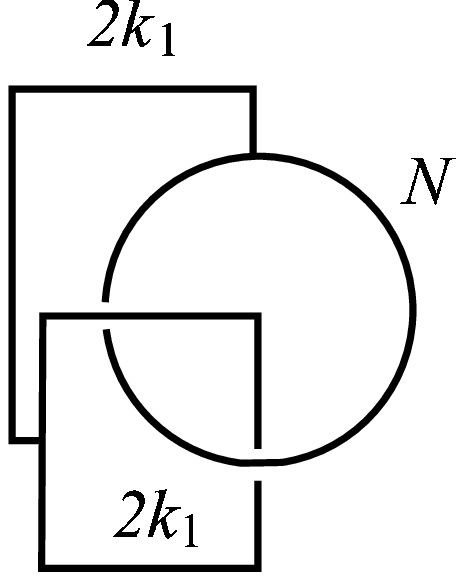}}
}$.\\

\noindent
Comparing the two expressions above, we obtain
$J_{14_{41721}}(A; N)=J_{14_{42125}}(A; N)$.
\qed

By using the same argument as for the previous pair
we showed that $14_{41763}$ and $14_{42021}$ have the same
CJP. This pair has again the same 2-cable HOMFLY polynomials,
but different Whitehead double HOMFLY polynomials.

\begin{figure}[h] 
\ig[width=4cm,height=4.5cm]{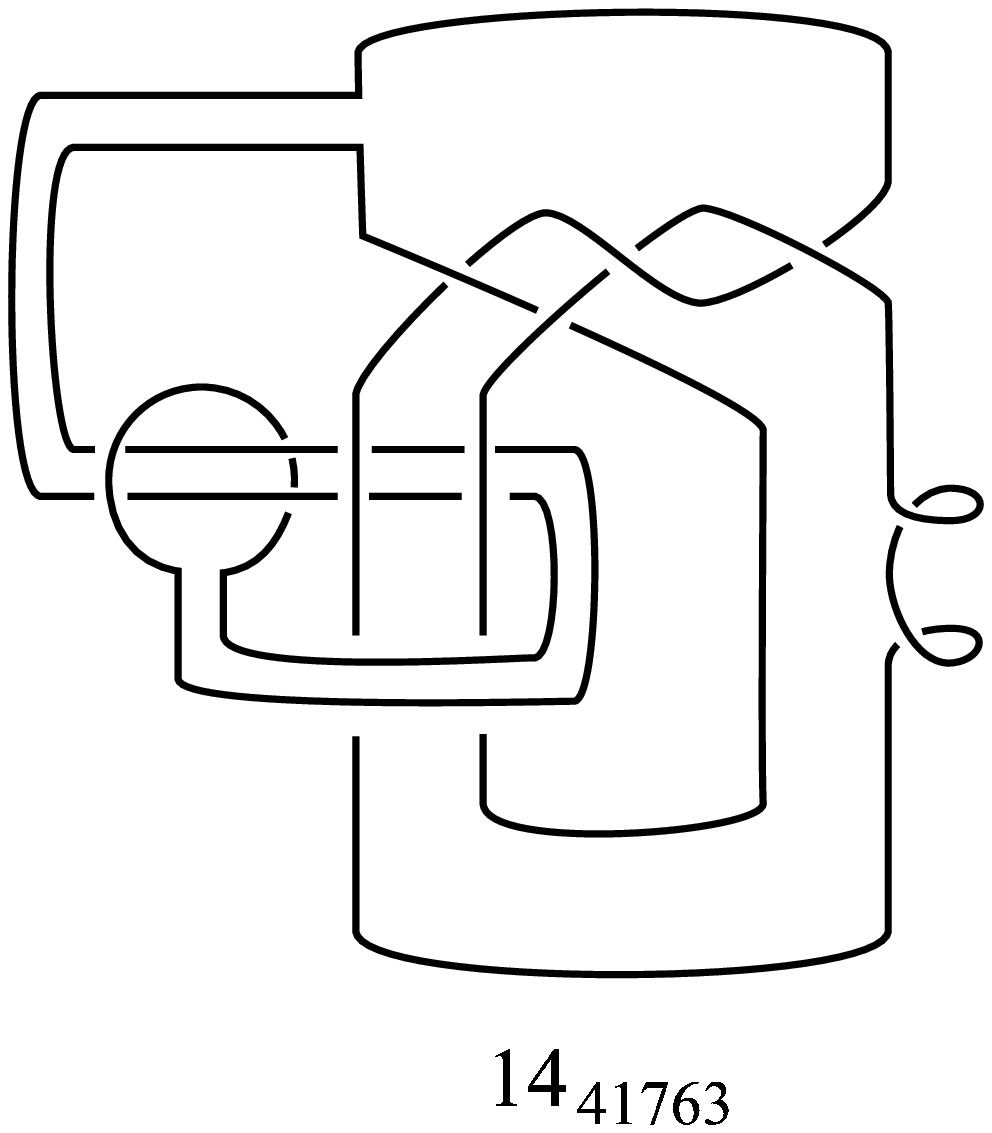}\hspace{5mm}
\ig[width=4cm,height=4.5cm]{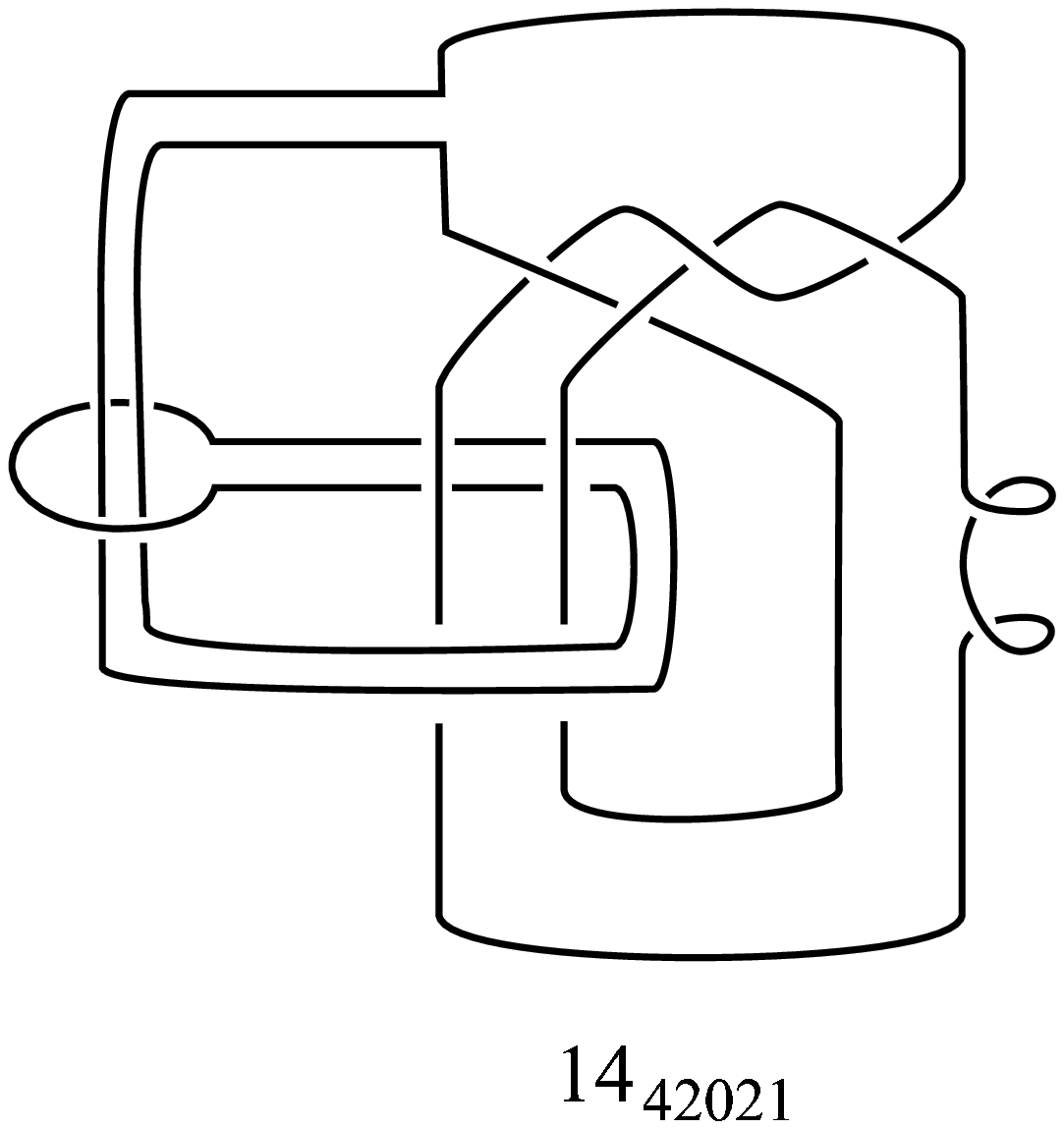}
	\caption{}
	\label{Pic16}
\end{figure}

For the proof of $J_{14_{41763}}=J_{14_{42021}}$, one needs
to find the suitable diagrams, shown in figure \reference{Pic16}.

Then we obtain

\noindent
$J_{14_{41763}}(A; N)$
\noindent
$=\displaystyle \sum_{k_{1}=0}^{N}\frac{<2k_{1}>}{<N,N,2k_{1}>}\frac{<2k_{1}>}{<2k_{1},2k_{1},2k_{1}>}
\br{\raisebox{-1.5cm}{\ig[width=3cm,height=2.5cm]{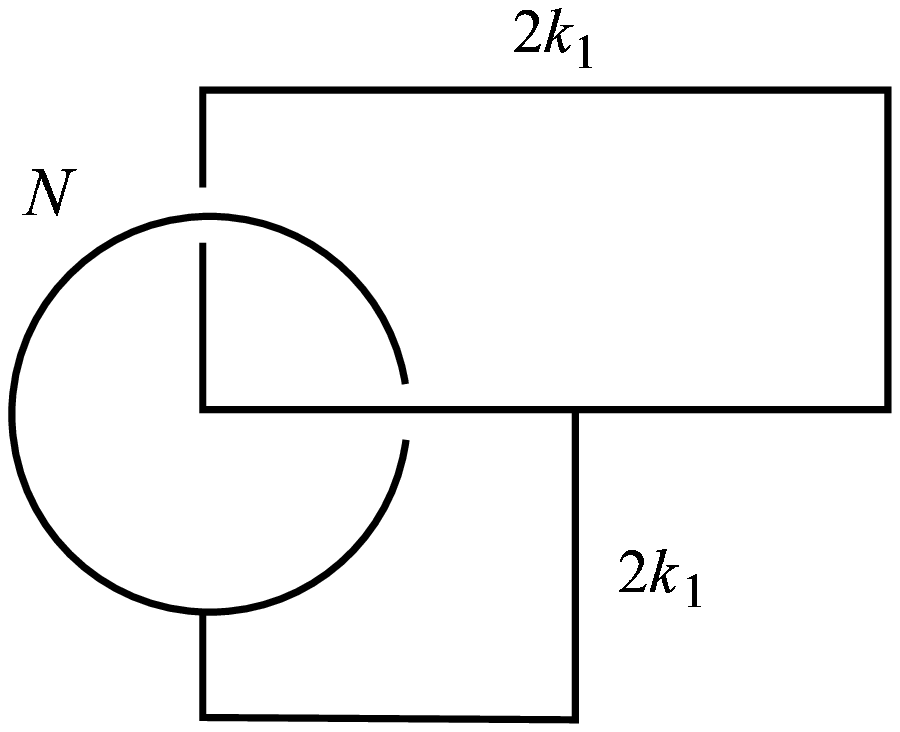}}}
\ %\times\\ \times
\br{\raisebox{-1.5cm}{\ig[width=2.5cm,height=3cm]{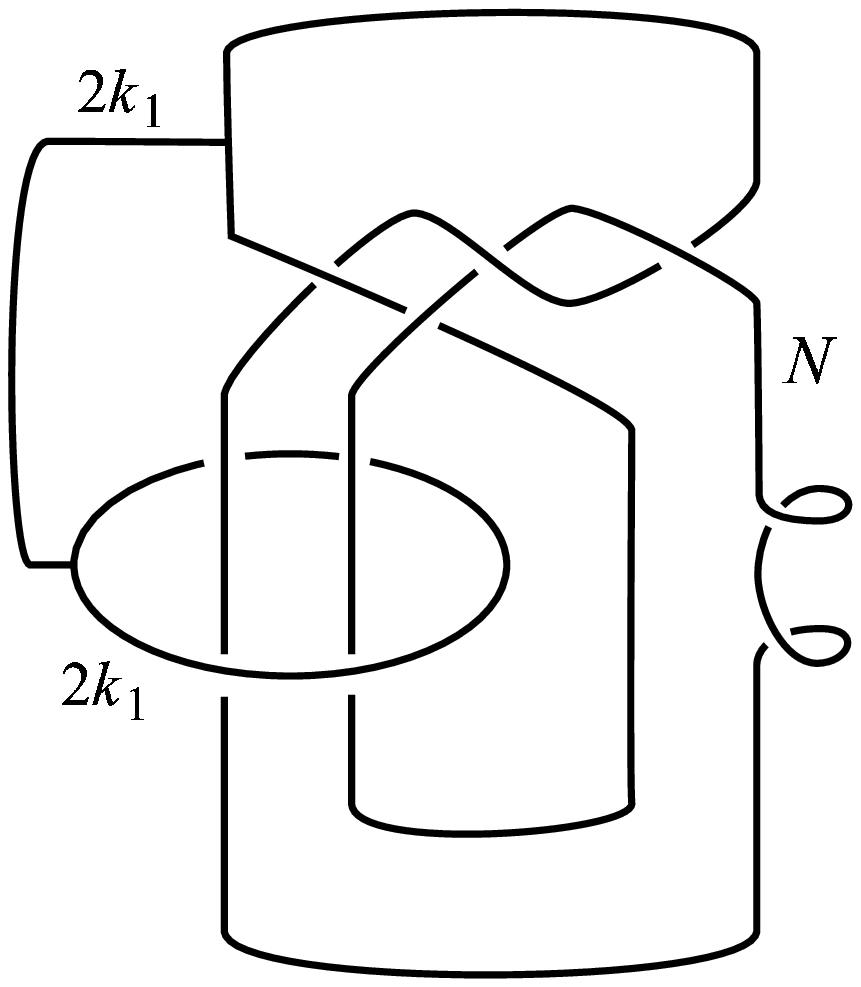}}}$.\\

\noindent
$J_{14_{42021}}(A; N)$
\noindent
$=\displaystyle \sum_{k_{1}=0}^{N}\frac{<2k_{1}>}{<N,N,2k_{1}>}\frac{<2k_{1}>}{<2k_{1},2k_{1},2k_{1}>}
\br{\raisebox{-1.5cm}{\ig[width=2.5cm,height=3cm]{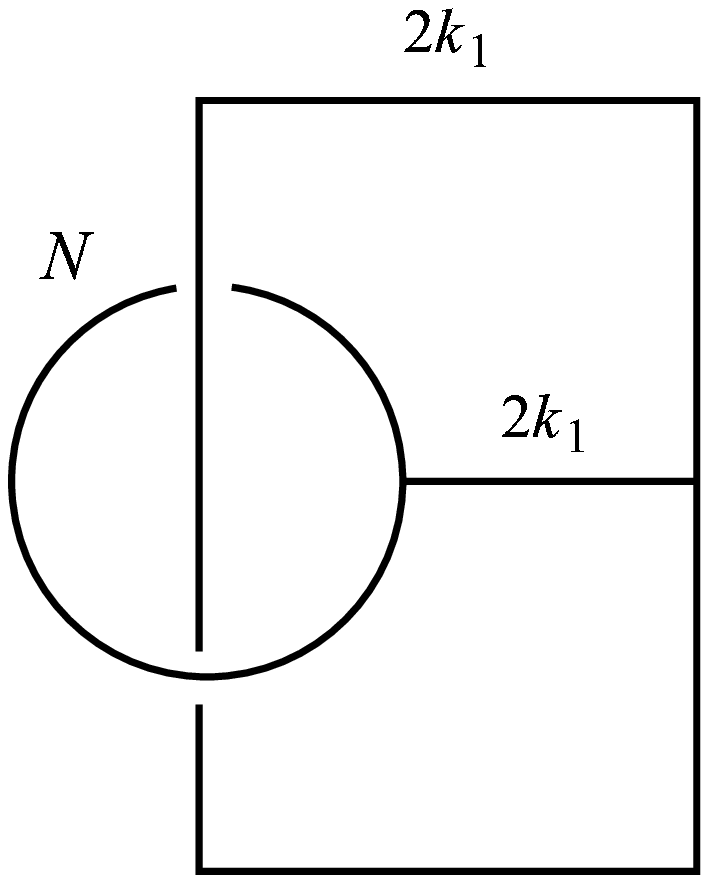}}}
\ %\times\\ \times
\br{\raisebox{-1.5cm}{\ig[width=2.5cm,height=3cm]{eps/Pic21.eps}}}$. \\

% \noindent
% Comparing two polynomials above we obtain $J_{14_{41763}}(A; N)=J_{14_{42021}}(A; N)$.\\

\noindent

Both $14_{41763}$ and $14_{42021}$ are concordant to the positive
trefoil knot. So we may put some weight on the following problem:

\begin{problem}\label{ph}
Does CJP detect a prime, simple knot modulo mutation and concordance?
More exactly, are two prime, simple knots with the same CJP
(a) mutants or concordant, or (b) related by a sequence of
concordances (preserving the CJP) and mutations?
\end{problem}

\begin{rem}\label{r44}
The invariance proof of the CJP exhibits a common pattern. Our
pairs are obtained by a certain involution of a 3-string tangle,
two of whose strands are parallel. In \cite{Gal} this kind of
operation was considered as a kind of what is called a ``cabled
mutation'', which in turn is
explained to be a special case of Ruberman's \cite{Ruberman}
(2,0)-mutation. Ruberman's work then shows why the knots in
theorem \ref{rg} have equal volume. It had been realized,
already in \cite[page 310, line -22]{CL}, that also the CJP is
invariant under a (2,0)-mutation (and hence cabled mutation).
Our examples turn out to be instances of this phenomenon, an
insight we owe to the explanation in \cite{Gal}. Contrarily, we
stress that at least our initial pair was quoted in \cite{Gal}.
% So it seems fair to say that we are the first to answer question
% ref{q1.1'}.
\end{rem}

\subsection{Knots with equal CJP, distinguished by
HOMFLY and Kauffman polynomials\label{S32}}

As a follow-up to proposition \reference{th12}, we make the
following observation:

\begin{rem}\label{r35}
Note that the colored Jones polynomial of $K$ determines the colored
Jones polynomial of the Whitehead doubles of $K$. (This is clear
from Kauffman bracket skein theory, but an exact formula is given
in \cite{T}.) So we see that the Whitehead doubles of our pair
are examples of knots with the same colored Jones polynomial but
different HOMFLY and Kauffman polynomial. (See remark \ref{rK},
and compare also part (1) in theorem \reference{th1.2}.)
\end{rem}

The only such knots we initially thought of are 3-satellites of
mutants. There have been some attempts to distinguish mutants
by applying polynomials on their 3-satellites. For the HOMFLY
polynomial the first such computation (and distinction) was
apparently done for 3-cables of the Kinoshita-Terasaka (K-T) and
the Conway (C) knot by Morton, at the time \cite{MorTra} was written.
Some further account is given in \cite{MorCro,MorRy}, incl. a
calculation using the quantum group (which is essentially equivalent
to evaluating parts of the HOMFLY polynomial). This work yields
different examples illustrating part (1) in theorem \ref{th1.2} for
the HOMFLY invariants. These examples address the 2-cable HOMFLY
polynomial, by taking the companions to be 2-cables of mutants.
A disadvantage of the way we collected our examples (see the
beginning of \S\reference{S5}) is that they pass by this case.
Contrarily, we deal with the Kauffman polynomial, and with
the double cover (which will be treated in our appendix).
For the Kauffman polynomial mutant 3-satellite calculations might
have been attempted, but despite some quest, we %the first author
found no reliable record available.

A further example for the HOMFLY polynomial is provided
in \cite{Gal} (the knot in Figure 3.5 therein). Still it is
worth noticing that the variance of the HOMFLY polynomial under
a cabled or (2,0)-mutation was understood already from Morton's
(aforegoing) above explained work, as credited also in \cite{CL}.
Satellite knots (in opposition to links) deserve no particular
prominence: by trivial skein theory standard cable polynomials are
linear combinations of connected cable polynomials. Therefore,
if former distinguish something, so must do (some of the) latter%
\footnote{This is something very different from distinguishing whether
the companion is a knot or link!}. On the other hand, our Whitehead
double pair from the proof of proposition \ref{th11} is somewhat
simpler than the knots in both \cite{Gal} and \cite{MorCro}.
In \cite{Gal} also similar examples for the Kauffman polynomial
were expected. The related Question 1.6 therein can be answered
(negatively, for the Kauffman polynomial), for instance,
with the 2-cable distinction of $14_{41721}$
and $14_{42125}$ in remark \ref{rK}.
%(where the polynomial can also be completely calculated).
After us, other pairs (distinguished by
either 2-variable polynomials) were given in \cite{MorRy2}.

We will give later several further, more subtle, examples that
relate also to the other mutant properties concerning 2-satellite
polynomials. As far as only Question \ref{q1.1'} is considered,
one can construct easily an infinite sequence of examples.

\proof [of Theorem \reference{rg}]
Figure \reference{Pic1} gives a ribbon presentation of the knots
as a pair of disks connected by a band.
We add an equal number of full-twists
in the bands. The proof that the CJP coincides is essentially the same,
only the half-twist coefficients of \cite{MV} enter additionally
into the formula. The property almost all knots to be simple knots is
established most easily using Thurston's hyperbolic surgery theorem
(see for example \cite{Thurston1,Thurston2}).
This theorem implies that all but finitely many of the knots will
be hyperbolic (and that when the number of twists goes to infinity,
the hyperbolic volumes will converge to the volume of the limit
link, which is $\approx 10.99$).

The distinction using the Whitehead double skein polynomial can
be argued about as follows. Since the polynomial distinguishes
our particular example, there is a Vassiliev invariant $v$ of some
degree $n$ contained in the Whitehead double skein polynomial (of
fixed framing), differing on both. (From Table \reference{figP},
one can calculate that $n=11$ is enough.) Now a Vassiliev
invariant is known to behave polynomially in the number of
twists. (This topic was initiated in \cite{Tr} and then expanded
\cite{Bseq}.) So $v$ will distinguish all but at most $n$ of
the pairs. Since $v$ is determined by the Whitehead double
HOMFLY polynomial, latter will differ too, and thus the pairs
are not mutants. \qed

\begin{rem}
An easy skein argument shows that the HOMFLY and Kauffman
polynomial coincide in all pairs. It will require a bit
more labor to check this for the 2-cable HOMFLY polynomial
(and we have not done so). %, while we do not know how to check
% it for the volume (though it is asserted by the Volume conjecture).
\end{rem}

% Using the multiplicativity of the CJP under connected sum, and
% either Ruberman's work, or a similar distinction argument based on
% Vassiliev invariants, we have

\begin{corr}\label{crn}
For any number $n$, there exists a family of $n$ distinct knots
with equal CJP, which are not mutants or satellites thereof. 
\end{corr}

\proof For given $n$, and with $K_1=14_{41721}$ and $K_2=14_{42125}$,
consider the knots $\{\,\#^k\,K_1\,\#\,\#^l\,K_2\,:\,k+l=n-1\,\}$.
By the multiplicativity of the CJP under connected sum, all these
knots have the same CJP. To prove that they are not mutants,
we use the following well-known fact: if a composite
knot has a mutant, then the mutant is also composite, and
prime factors of eithers correspond $1-1$ up to mutation. \qed

% \begin{rem}\label{rty}
It is intuitively clear that coincidences of the CJP are rather
sporadic. Thurston's results that the hyperbolic volumes of
3-manifolds form a well-ordered set and are finite-to-one
invariants imply that only finitely many knots with no Seifert
fibered pieces in the JSJ decomposition have the same Gromov norm.
The Volume conjecture then asserts that only finitely many such
knots have the same CJP, and suggestively the same is true under
dropping the mild restriction on the JSJ decomposition. So the
examples in the corollary are likely the most one could obtain.
% See also remark \reference{rpr}.
% \end{rem}

% \begin{rem}\label{red}
% The insight in remark \ref{rty} is highly non-constructive, and 
Contrarily, we have no procedure to find all the knots for given CJP.
It is known, though, that there are non-concordant mutants (see
for example \cite{KL}), so, like mutation, concordance alone
will not suffice to relate knots with the same CJP (see problem
\reference{ph}).
% \end{rem}

\section{Mirror images\label{S4}}

Since in the tables of \cite{HT} a knot $K$ is considered
equivalent to its \em{mirror image} (or \em{obverse}) $!K$,
a careful detection of mutants entails also attention to
chiral knots which are mutants to their obverses. An example
of a 16 crossing knot, given by Sakuma and Kanenobu, is quoted
in \cite{Kauffman}. No such knots occur up to 13 crossings,
because the 2-cable HOMFLY polynomial distinguishes all
chiral knots from their obverses, although for some knots
the uncabled HOMFLY and Kauffman polynomial fail, most prominently
$9_{42}$, and also $10_{71}$. Note that taking the mirror image
replaces one of the variables in the HOMFLY, Kauffman or
Colored Jones polynomial by its inverse. Polynomials that
remain invariant under this interchange are called below
\em{reciprocal} or \em{self-conjugate}.

For 14 crossings, the 2-cable HOMFLY polynomial fails distinguishing
mirror images on 15 prime chiral knots on which the uncabled polynomials
fail too. 13 of them are indeed found to be mutants to their obverses,
because they are mutants to achiral knots. The other two knots are
distinguished by the Whitehead double HOMFLY polynomials. These knots
are $14_{3802}$ and $14_{29709}$; see figure \reference{Picach}.

\twofigures{
\begin{array}{c@{\hspace{8mm}}c}
\ig[width=4.0cm,height=3.5cm]{eps/t1-14_3802.eps} &
\ig[width=4.0cm,height=3.5cm]{eps/t1-14_29709.eps} \\[4mm]
14_{3802} & 14_{29709}
\end{array}
}{
\begin{array}{c}
\ig[width=4.0cm,height=3.5cm]{eps/t1-staple2.eps} \\[4mm]
S
\end{array}
}
{}{The staple tangle}{Picach}{Picstaple_}

It is natural to ask how well the Colored Jones polynomial would
distinguish mirror images of such difficult cases. So we considered
the knots we did not find mutants to their obverses, but which
have zero signature and self-conjugate uncabled polynomials.
Note that the Kauffman polynomial determines for \em{knots} the
2-Colored Jones polynomial by a result of Yamada \cite{Yamada}
(see also \cite{King}).
So for knots with reciprocal uncabled polynomials, the 1- and
2-Colored Jones polynomials must be reciprocal too. We tried to
determine the 3-colored polynomial using the \verb|KnotTheory`|
Mathematica Package of Dror Bar-Natan \cite{BN2} (and its Colored
Jones polynomial facility co-written with S.~Garoufalidis). We
obtained the polynomial only for a handful of our 14 crossing knots
(and no knots of more crossings).

\begin{exam}\label{x41}
Among others we calculated that the 3-Colored Jones polynomial
distinguishes $14_{29709}$ and its mirror image. This
polynomial is not reciprocal, thus providing the example for
part (2) in theorem \ref{th1.2}:

%% [INCLUDE CALCULATION FOR $14_{29709}$  \& MIRROR IMAGE ]
 $J_{14_{29709}}(3;q) =
 1/q^{31} - 2/q^{30} - 1/q^{29} + 8/q^{27} + 5/q^{26} - 18/q^{25} -
 21/q^{24} + 16/q^{23} + 64/q^{22} + 3/q^{21} - 108/q^{20} - 
 76/q^{19} + 140/q^{18} + 194/q^{17} - 105/q^{16} - 353/q^{15} - 
 5/q^{14} + 483/q^{13} + 217/q^{12} - 569/q^{11} - 468/q^{10} + 
 560/q^9 + 734/q^8 - 480/q^7 - 957/q^6 + 346/q^5 + 
 1116/q^4 - 187/q^3 - 1208/q^2 + 24/q + 1240 + 132 q - 1208 q^2 - 
 290 q^3 + 1119 q^4 + 442 q^5 - 967 q^6 - 571 q^7 + 
 754 q^8 + 649 q^9 - 493 q^{10} - 661 q^{11} + 240 q^{12} + 
 575 q^{13} - 17 q^{14} - 437 q^{15} - 110 q^{16} + 262 q^{17} + 
 158 q^{18} - 121 q^{19} - 129 q^{20} + 24 q^{21} + 82 q^{22} + 
 9 q^{23} - 31 q^{24} - 17 q^{25} + 8 q^{26} + 9 q^{27} - q^{28} - 
 q^{29} - 2 q^{30} + q^{31}$\,.
\end{exam}

Similarly, the 3-colored polynomial
distinguished the few other knots we could evaluate
it on from their mirror images. Contrarily, we can prove now the
following, which also settles part (3) (b) of theorem \ref{th1.2}.

\begin{thm}\label{thach}
There exist infinitely many hyperbolic knots $K_n$, such that
$K_n$ has the same CJP as its mirror image, but they are not mutants.
\end{thm}

\proof %[of Theorem \reference{thach}]
Let $S$ be the (3-string) staple tangle in figure \ref{Picstaple_}.
(The name is taken from \cite{SW}, where similar tangles
are heavily used.) If one turns around one of the diagrams
in Figure \reference{Pic1}, one sees that the knots differ by a
staple turn. This is a rotation of $S$ by $\pi$ around the axis
vertical in the projection plane. It has the effect of mirroring
the staple. The calculation for the proof of proposition
\ref{th12} shows that a staple turn does not change the CJP.

Now consider the knot $K_n$ on the left in figure \ref{Picstaple}. The
diagram contains a pretzel tangle $T_n=(-n,n)$ with $n$ odd (we showed
it for $n=3$), and the remaining part consists of the join of two
staple tangles $S$. For us it is useful here that
mutation of $T_n$ along the axis horizontal in the projection plane
turns it into its mirror image. So $K_n$ can be transformed into
$!K_n$ by two turns of tangles $S$ and one mutation (of $T_n$).
Thus $K_n$ and $!K_n$ have the same CJP. The limit link $L_\infty$,
obtained by placing circles $N_{1,2}$ around each of the groups of
$n$ twists (and then be ignoring the twists by \cite{Adams}),
is a hyperbolic link (of volume $\approx 28.07$). Then Thurston's
hyperbolic surgery theorem assures that $K_n$ are hyperbolic for
large $n$.

\begin{figure}[h] 
\[
\def\gr#1#2{
  \diag{1cm}{4.6}{4.2}{
    \picputtext{2.3 2.1}{
       \hbox{\ig[width=4.6cm,height=4.2cm]{eps/#1.eps}}
    }
    #2
  }
}
\begin{array}{c@{\hspace{9mm}}c@{\hspace{9mm}}c}
\gr{t1-staple3_3}{} &
\gr{t1-staple3}{} &
\gr{t1-staple3_lim}{
  \picputtext{1.01 2.0}{$N_1$}
  \picputtext{1.11 0.4}{$N_2$}
}
\\[22mm]
K_3 & K=K_{\pm 1} & L_\infty
\end{array}
\]
\caption{} \label{Picstaple}
\end{figure}

It remains to prove that $K_n$ is not mutant to $!K_n$. Let $L'$ be the
2-component trivial sublink of the limit link $L_\infty$, consisting
of the two circles $N_{1,2}$. Then $\pm 1/m$ surgery along the
components of $L'$ performs at two occasions $m$ full-twists
inside the tangle $T_n$, turning it into $T_{2m+n}$. The
work in \cite{Bseq} then implies that for a Vassiliev invariant
$x$, the function $n\mapsto x(K_n)$ is a polynomial function in
$n$, of degree at most $d:=\deg x$. Now (quite non-trivial)
calculation of the Whitehead double HOMFLY polynomial of $K:=
K_{\pm 1}$ found a Vassiliev invariant $x$ with $x(K)\ne x(!K)$.
Thus $x$ will distinguish from their mirror images all but
at most $d$ of the $K_n$. Since $x$ is not changed under mutation,
we see that $K_n$ and $!K_n$ are not mutants. \qed

\begin{rem}\label{ry}
It is known that mutation preserves the double branched cover $M_2(K)$
of a knot $K$, and this fact will become very important later.
A main motivation for the construction of the examples here was to
manifest the usefulness of polynomial and Vassiliev invariants
as a tool to exclude mutation, in opposition to the study of $M_2(K)$.
It is clear that $M_2(K)=M_2(!K)$, and since the peripheral system
of $K$ loses orientation in $M_2(K)$, there seems no way to prohibit
a mutation between $K$ and $!K$ by studying $\pi_1(M_2(K))$ (alone).
The other type of invariants coming to mind are those derived from
the Blanchfield pairing on the Alexander module (signatures and
linking forms). But one easily comes across knots of trivial Alexander
polynomial, where such methods will fail either. (See e.g. the knot
$K$ in figure \reference{Picstaple}.)
\end{rem}

\begin{exam}
The knot $14_{3802}$ in figure \reference{Picach} is
% an instance that creates difficulty on the opposite side to $K$ of
% figure \reference{Picstaple}. It is 
particularly interesting,
because it is alternating. For alternating knots one can apply
the strong geometric work in \cite{MenThis} and \cite{Menasco}.
Former result shows that $14_{3802}$ is chiral, and using the
latter result one can easily deduce that $14_{3802}$ has no mutants,
so in particular it cannot be a mutant to its mirror image. This
conclusion is here especially difficult to obtain using the
polynomials. We could not decide (calculationally) whether
$14_{3802}$ and its mirror image have the same Colored Jones
polynomial (or even the same 3-Colored Jones polynomial).
\end{exam}

\begin{exam}
$14_{29709}$ and its mirror image provide an example of knots
with the same 2-cable HOMFLY polynomials and the same volume,
but different 3-Colored Jones polynomials. For knots of different
volume, but equal 2-cable HOMFLY polynomials, we found previously
the pair ($12_{341}$, $12_{627}$) in \cite{cab}, and checked now
similarly that the 3-Colored Jones polynomials are different.
These examples show that the Colored Jones polynomial is not
determined by the 2-cable HOMFLY polynomial.
\end{exam}

In contrast, the following question remains open:

\begin{question}
Are there examples of knots with different (3-)colored Jones
polynomial but equal Whitehead double HOMFLY (or Kauffman) polynomials?
\end{question}

In \S\reference{S3} we observed that the Colored Jones polynomial
does not determine in turn even the uncabled HOMFLY polynomial.

Another interesting question raised by our verification is:

\begin{question}
Is every knot which is a mutant to its obverse actually
a mutant to an achiral knot?
\end{question}

The above explanation implies that it is so up to 14 crossings.
Similarly, there are 288 chiral prime knots of
16 crossings which the (uncabled) HOMFLY and Kauffman polynomials
and the signature fail to separate from their mirror image. 117 are
mutants to achiral 16 crossing knots, and all the other 171 are 
ruled out by the Whitehead double HOMFLY polynomials (though the 
2-cable HOMFLY polynomials fail on 6 of them, all non-alternating).
For 15 crossings all chiral prime knots are found distinguished from
their obverses by one of the uncabled polynomials or the signature. 

\section{More difficult examples\label{S5}}

The pairs we presented came up in the first author's project to
determine mutations among the low crossing knots tabulated in 
\cite{HT}. Up to 13 crossings this task was completed by tracking 
down coincidences of Alexander, Jones polynomial and volume (up to a
certain computable precision) on the one hand, and then exhibiting
the mutation in minimal crossing diagrams on the other hand.

In contrast, a (non-complete) verification of 14 and 15 crossing 
knots exhibited several more difficult cases, discussed in \cite
{mut}. For some pairs we had to seek non-minimal crossing diagrams 
to display the mutation. Others provided examples of the type we 
showed in \S\reference{S3}. Here may be a proper place to stress 
that the match of the $P$, $F$ and 2-cable $P$ polynomials on all 
our examples is intentional. Beside a mutation status check, it
served as a selector for good candidates with the same CJP. The proof 
to confirm equality of the CJP for each pair entails the quest for 
proper presentation of the knots, and so requires more effort than 
its length may let appear. We could therefore examine only a limited 
number of pairs. It is difficult also to exclude pairs, since a direct 
calculation of the 3-colored polynomial is feasible only in the
fewest cases.

Among the remaining, most problematic, pairs is $(14_{41739},14_
{42126})$, and a number of pairs of 15 crossing knots. An extensive
(though not exhaustive) search of diagrams up to 18 crossings
failed to show a mutation, but, along with all the invariants
that do so for the previous pairs, we established (in the way
explained before the proof of Proposition \ref{th11}) that Whitehead
double skein polynomials also coincide. 

By a similar calculation to the proof of Proposition \ref{th12}, we
managed to verify for some pairs that the colored Jones polynomials
are also equal. These pairs include the 14 crossing knots 
$14_{41739}$ and $14_{42126}$ in Figure \ref{Pic26}.

\twofigures{
\ig[width=3.5cm,height=3.9cm]{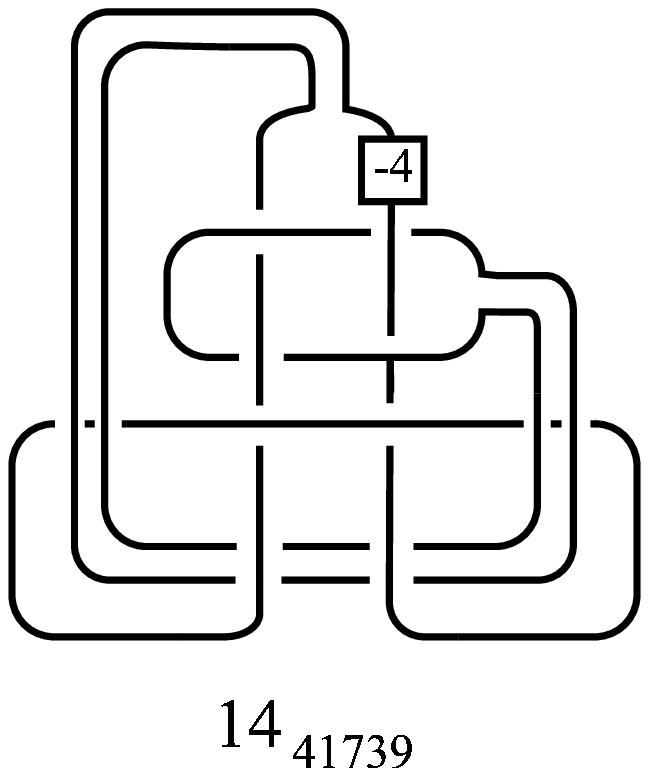}\hspace{5mm}
\ig[width=3.5cm,height=3.9cm]{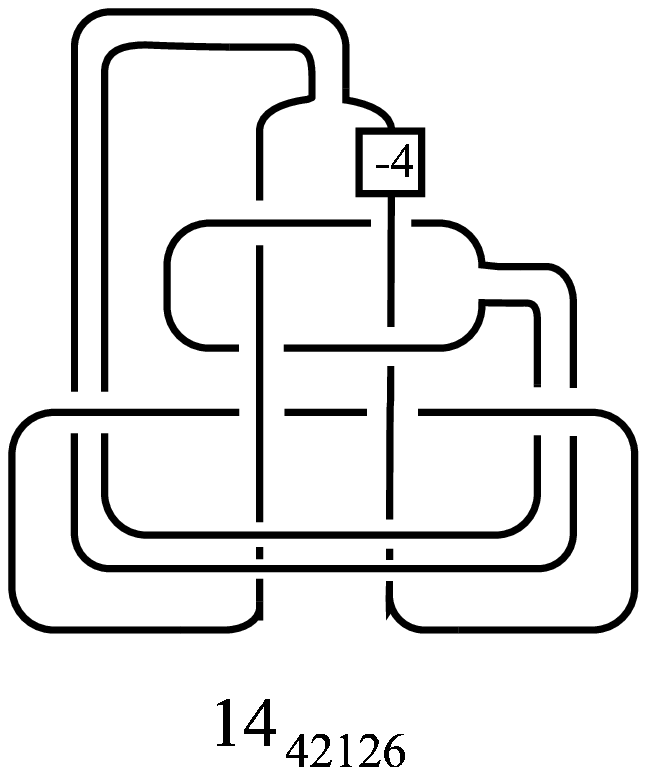}
}{
\ig[width=3.5cm,height=3.9cm]{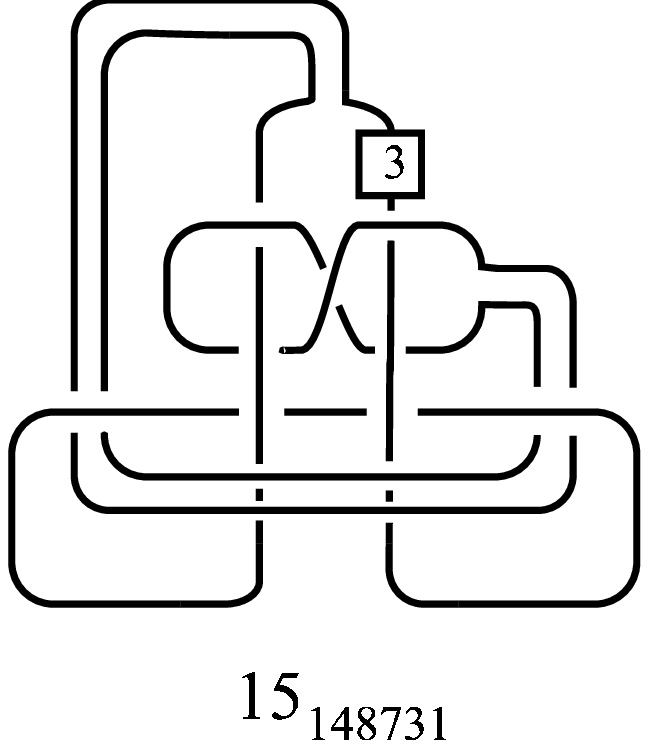}\hspace{5mm}
\ig[width=3.5cm,height=3.9cm]{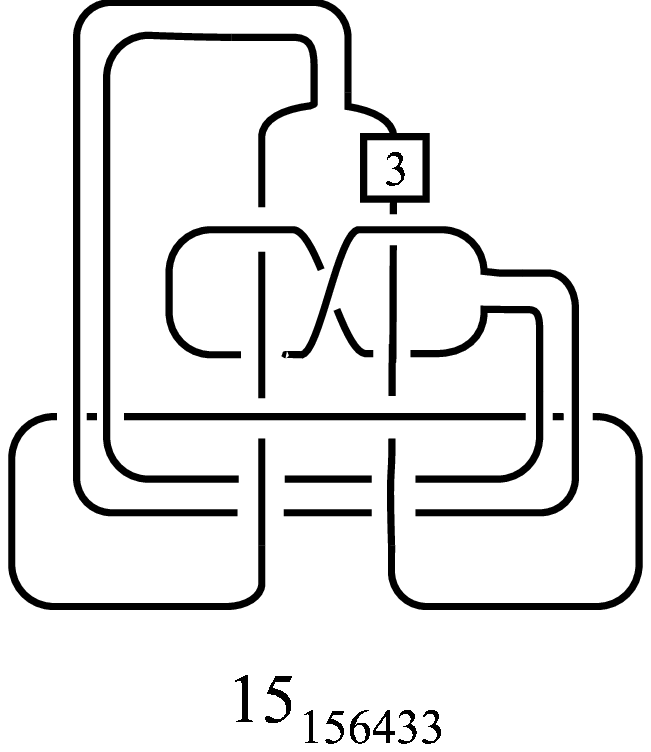}
}{}{}{Pic26}{Pic31}

\noindent
{\it Proof of $J_{14_{41739}}=J_{14_{42126}}$}.
We use the diagrams in the figure (a box with an
integer $w$ inside means $|w|$ kinks of writhe $\sgn(w)$). Setting
\begin{equation}\label{Gm}
\Gm(N,k_1,k_2)\,:= \frac{<2k_1>}{<N,N,2k_1>}
\frac{<2k_{2}>}{<N,N,2k_{2}>}\frac{<2k_1>}{<2k_1,2k_{2},2k_1>}\,,
\end{equation}
we find
\[
J_{14_{41739}}(A; N)
= \sum_{k_{1}=0}^{N} \sum_{k_{2}=0}^{N}\Gm(N,k_1,k_2) \times
\br{\raisebox{-1.5cm}{\ig[width=2.5cm,height=3cm]{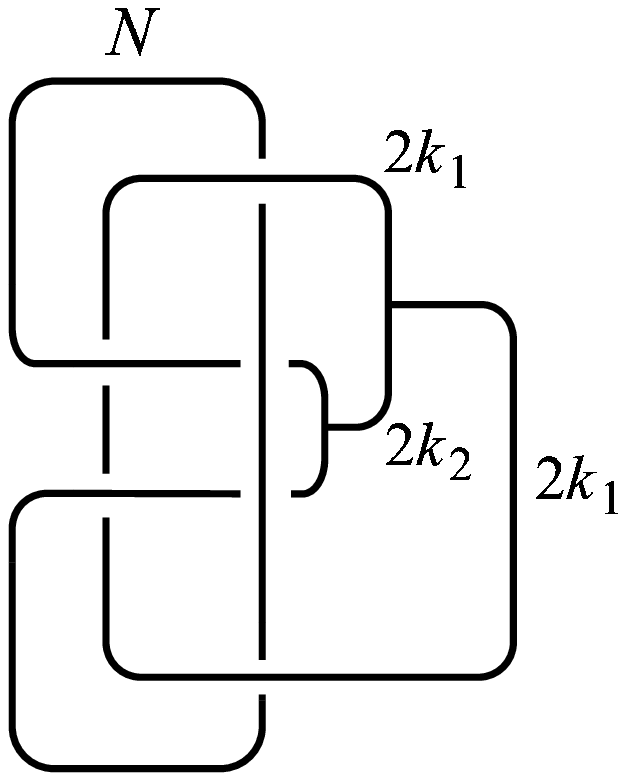}}}
\br{\raisebox{-1.5cm}{\ig[width=4cm,height=3cm]{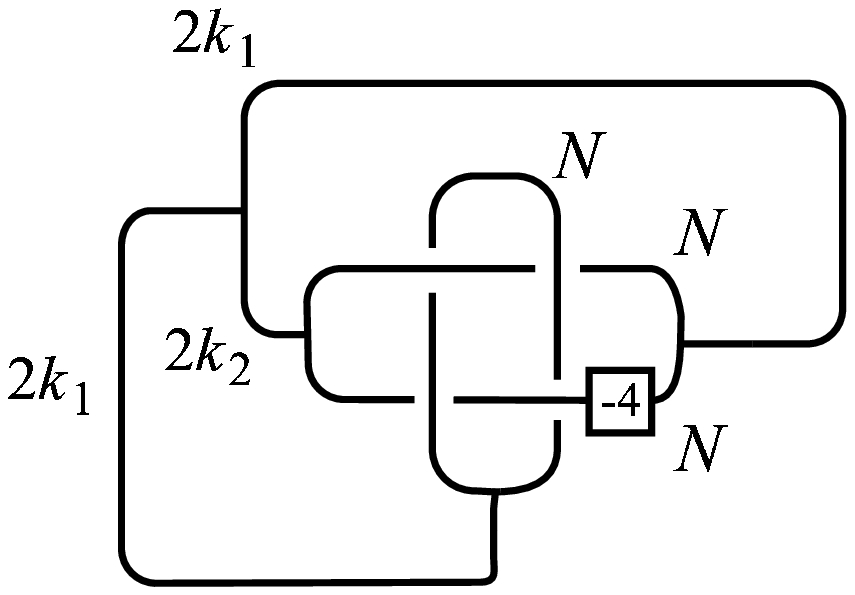}}}\,,
\]
\[
J_{14_{42126}}(A; N) =
\sum_{k_{1}=0}^{N} \sum_{k_{2}=0}^{N}\Gm(N,k_1,k_2) \times
\br{\raisebox{-1.5cm}{\ig[width=2.5cm,height=3cm]{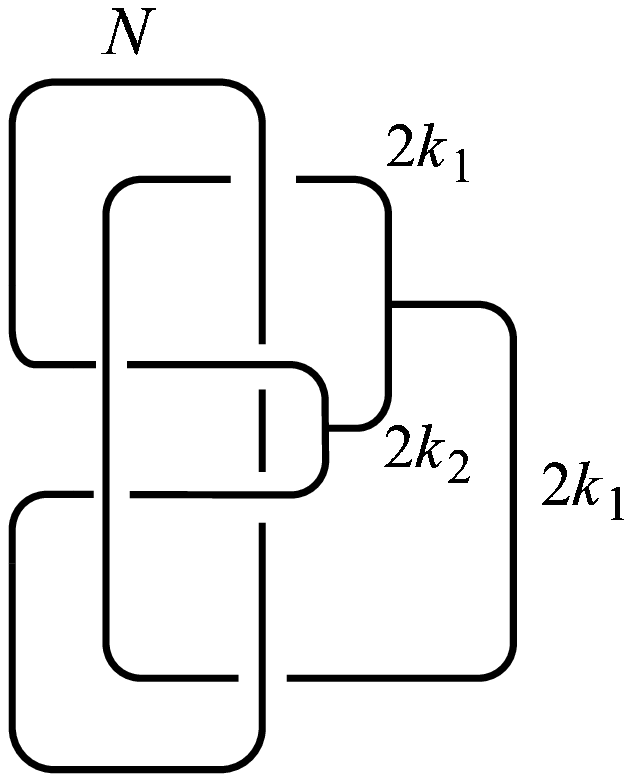}}}
\br{\raisebox{-1.5cm}{\ig[width=4cm,height=3cm]{eps/Pic29.eps}}}\,
\]\\[3mm]
(again with the sums restricted over $k_1,k_2$ for which
the denominator of $\Gm(N,k_1,k_2)$ makes sense). Comparison
shows again $J_{14_{41739}}(A; N)=J_{14_{42126}}(A; N)$. \qed

The CJP of the pair of (slice, trivial Alexander polynomial) knots
in Figure \ref{Pic31} is found equal from the displayed diagrams,
in a similar way to $14_{41739}$ and $14_{42126}$.

% {\it Proof of $J_{15_{148731}}=J_{15_{156433}}$}.

\begin{figure}[h]
\ig[width=3.0cm,height=5.5cm]{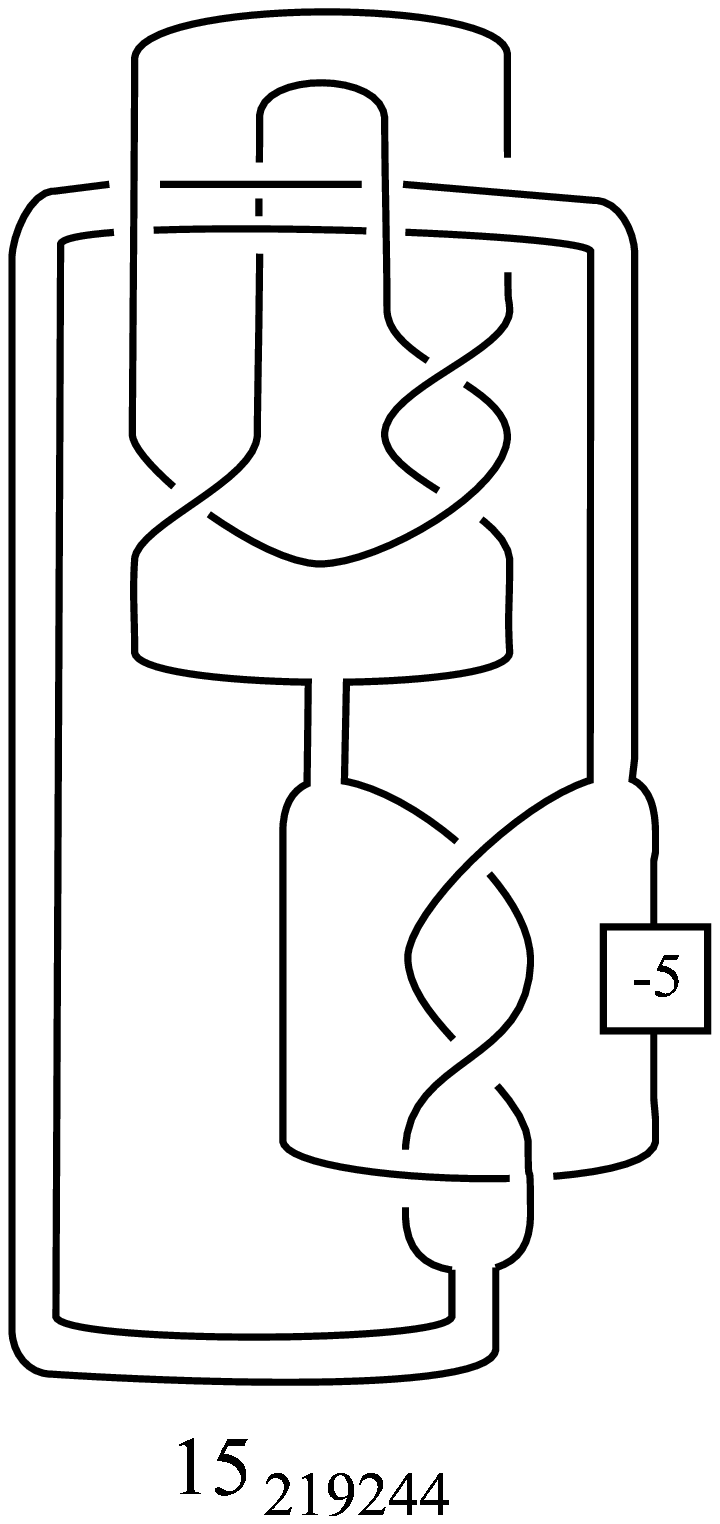}\hspace{5mm}
\ig[width=5.9cm,height=5.5cm]{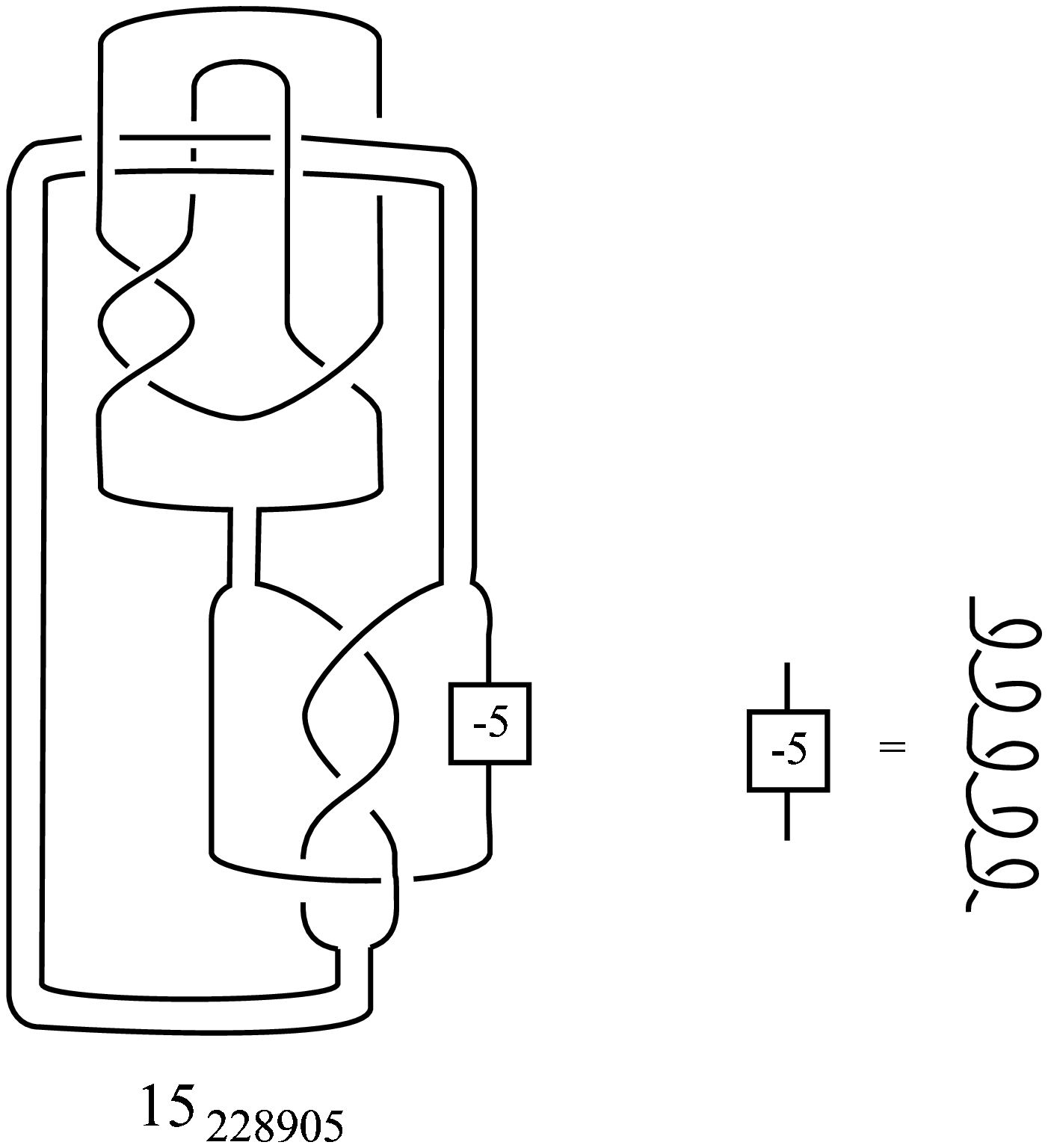}
	\caption{}
	\label{Pic18}
\end{figure}

A third pair is provided by the knots $15_{219244}$ and $15_{228905}$
in figure \reference{Pic18}.

\noindent {\it Proof of $J_{15_{219244}}=J_{15_{228905}}$}. With
$\Gm(N,k_1,k_2)$ as in \eqref{Gm}, we have from figure \ref{Pic18}
\[
J_{15_{219244}}(A; N)
= \sum_{k_{1}=0}^{N} \sum_{k_{2}=0}^{N}\Gm(N,k_1,k_2) \times
\br{\raisebox{-1.5cm}{\ig[width=2.5cm,height=3cm]{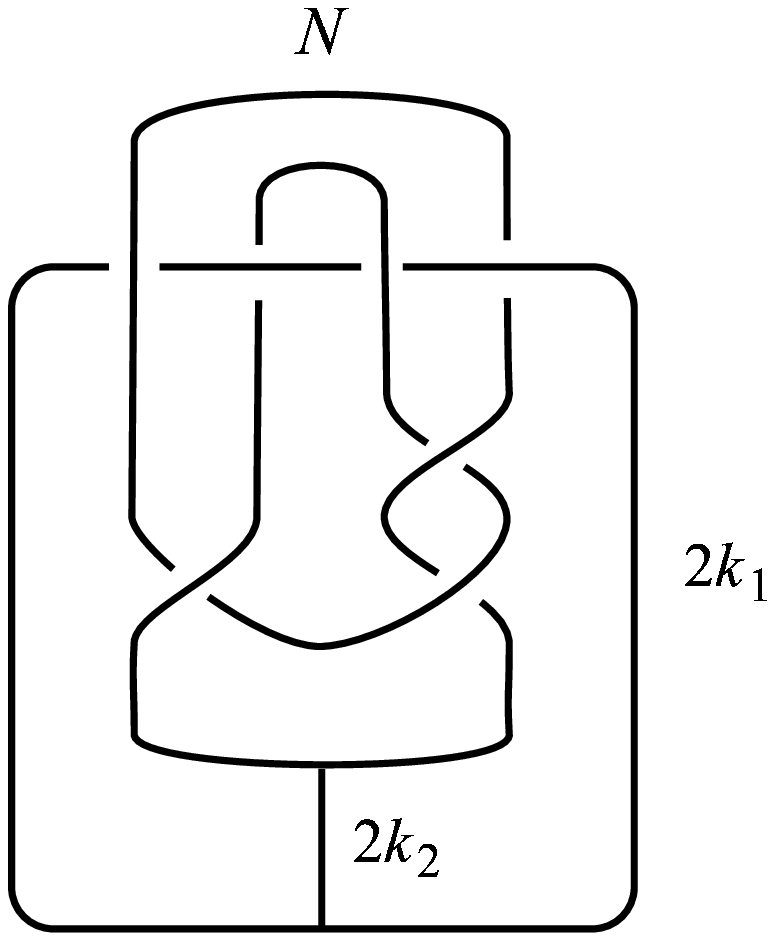}}}
\br{\raisebox{-1.5cm}{\ig[width=2.5cm,height=3cm]{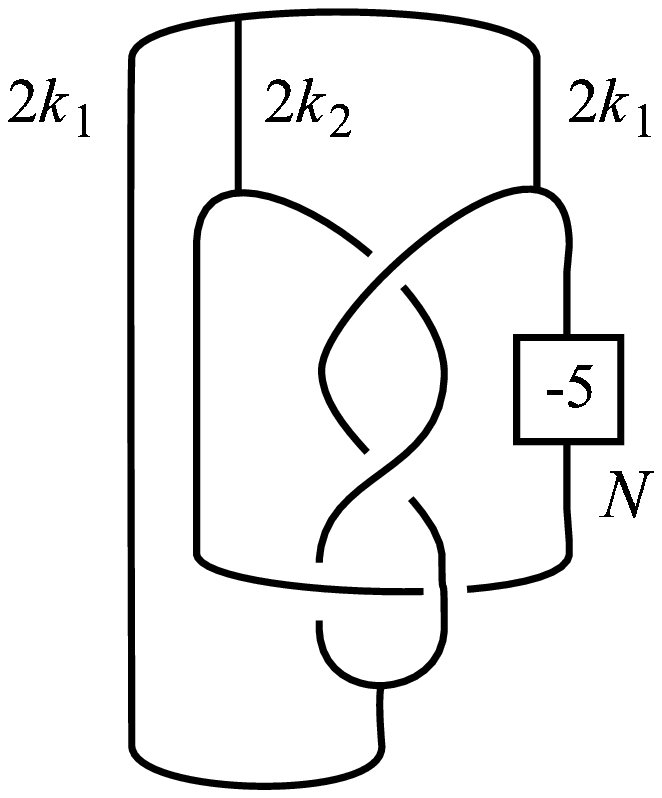}}}\,,
\]
\begin{equation}
% \tm
\def\tagform@#1{\maketag@@@ {\ignorespaces #1\unskip \@@italiccorr}}
\tag{\hbox to \z@{\hbox to \textwidth{\hss \qed}\hss}\phantom{\qed}} 
J_{15_{228905}}(A; N)
= \sum_{k_{1}=0}^{N} \sum_{k_{2}=0}^{N}\Gm(N,k_1,k_2) \times
\br{\raisebox{-1.5cm}{\ig[width=2.5cm,height=3cm]{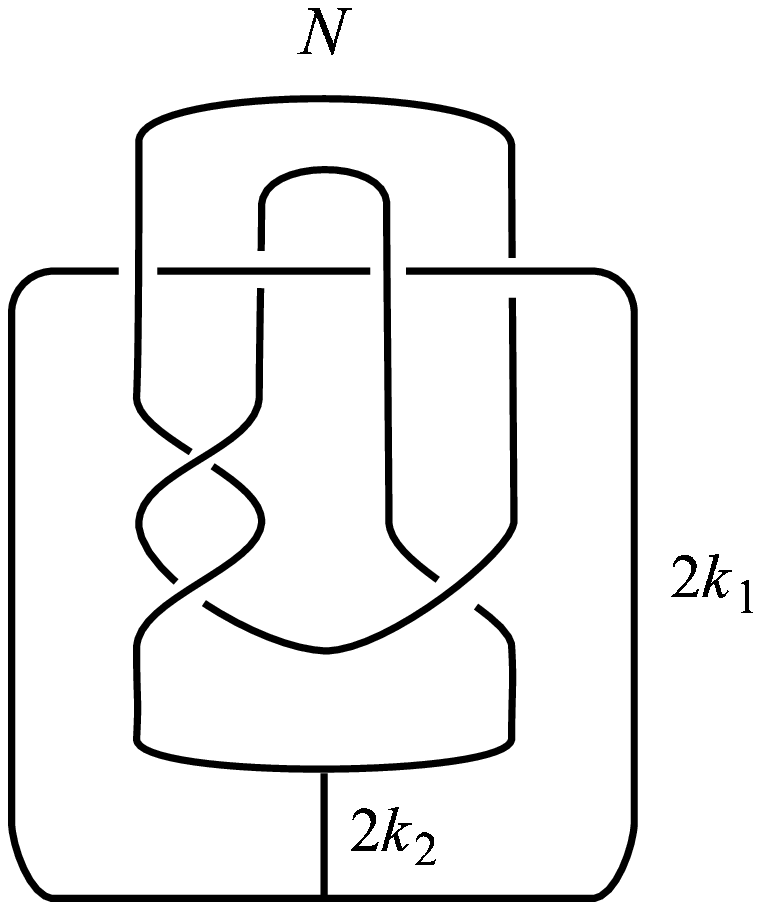}}}
\br{\raisebox{-1.5cm}{\ig[width=2.5cm,height=3cm]{eps/Pic24.eps}}}\,.
\end{equation}
% \noindent
% Comparing two polynomials above we obtain $J_{15_{219244}}(A; N)=J_{15_{228905}}(A; N)$.
% \qed
% [INCLUDE FEW DETAIL ON CALCULATION]

\begin{figure}[h] 
\[
\begin{array}{c@{\hspace{2cm}}c}
\ig[width=4.5cm,height=4.1cm]{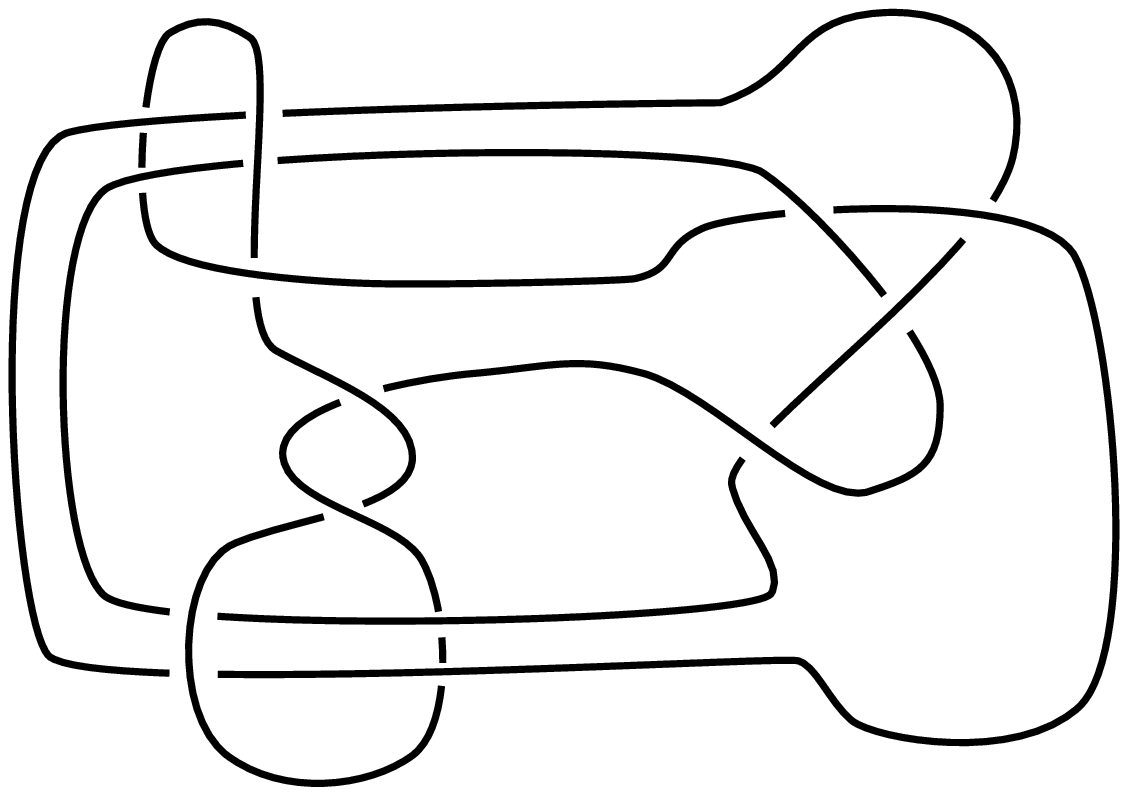} & 
\ig[width=4.5cm,height=4.1cm]{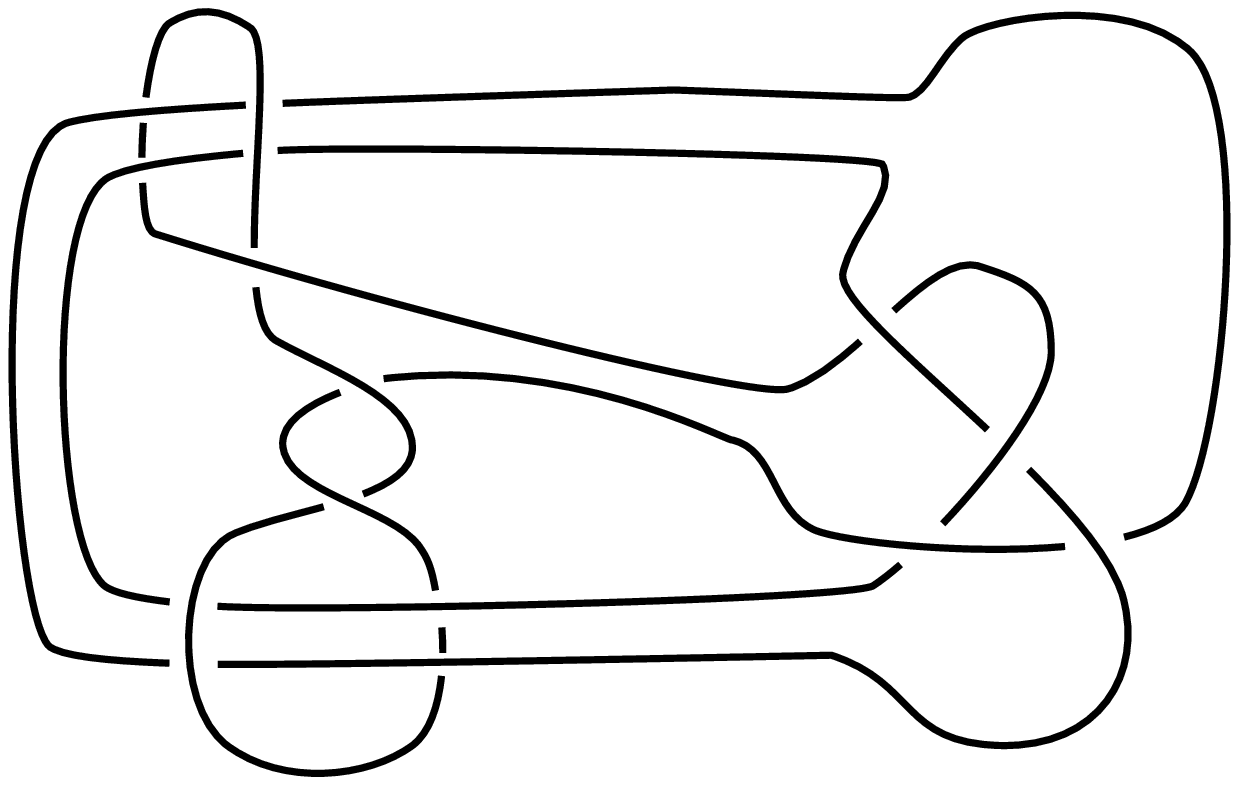} \\[-5mm]
15_{220504} & 15_{234873}
\end{array}
\]
\caption{}\label{Pic32}
\end{figure}

One more pair is $(15_{220504},15_{234873})$ in figure \ref{Pic32}.

The last two pairs shown are different from the previous two also
in the point that we do not know if they are concordant or not (see
problem \reference{ph}).
% (Note, however, that even if not, at least part (b) in problem
% \ref{ph} remains open, because of remark \ref{red}.)

In an attempt to exclude mutation, we relied on the fact that
mutants have the same double branched covers. We consulted
Daniel Matei, who succeeded in distinguishing these pairs
(and also the few others he tried, with differing degree of
effort) by the representations of the fundamental group of
the cover. In an appendix of the paper he gives some details
about his calculation.

The complexity of the 2-cable Kauffman polynomial makes its
evaluation very difficult. It is certainly not to be considered
a reasonable mutation criterion. We tried evaluating it on some
pairs (as in remark \ref{rK}), mainly driven by curiosity if
it would distinguish pairs left undistinguished by the HOMFLY
polynomials. We succeeded to determine the polynomial only for
two of the 15 crossing pairs with equal Whitehead double skein
polynomial. The 2-cable Kauffman polynomials, too, failed to
distinguish the knots.

\begin{exam}\label{x56}
One of these two pairs is $(15_{148731}$, $15_{156433})$ of 
Figure \reference{Pic31}. For this pair we were able to
calculate (and found to be equal) connected 2-cable Kauffman
polynomials for both mirror images, and thus we know that these
two knots satisfy \em{all} polynomial coincidence properties
known for mutants (those summarized in the introduction). The 
apparent lack of diagrams exhibiting the mutation deepened 
the decision problem whether the knots are mutants or not. 
Daniel Matei's subsequent exclusion result thus leads to
the most striking among the examples we have, showing complete
failure of the polynomial invariants to determine the mutation
status (see part (3) (a) in theorem \ref{th1.2}). On the
other hand, with theorem \reference{thach} we
saw examples, where the polynomial invariants seem indispensable
(see remark \reference{ry}). This underscores the significance
of both approaches.
\end{exam}

\begin{rem}
Although $\pi_1(M_2(K))$ becomes a major distiction tool, we know
of no pair in which its abelianizations $H_1(M_2(K))$ are different.
This relates to a question of Rong \cite[problem 1.91(4b)]{Kirby},
whether $J_K$ would always determine $H_1(M_2(K))$. This question
can be seen, together with the aforementioned relation of CJP to the
Alexander polynomial and the signature function, as part of a larger
conception whether the CJP might determine the Seifert matrix (up to
$S$-equivalence). For some related motivation see the remark after
problem 1.87 in \cite{Kirby}.
\end{rem}

\noindent{\bf Acknowledgement.}
The first author was supported by Postdoc grant P04300 of the Japan
Society for the Promotion of Science (JSPS). He would wish to thank
to his host Prof. T.~Kohno, and also to F.~Nagasato, H.~Morton
and P.~Traczyk for some helpful remarks. The second author's
research was supported by the 21st century COE program at the
Graduate School of Mathematical Sciences, the University of Tokyo
and Osaka City University. He would like to thank Prof. A. Kawauchi
and Prof. T. Kohno for helpful suggestions.

\begin{appendix}

% \renewcommand{\thesection}{\Alph{section}}
% {\def\@seccntformat#1{Appendix \csname the#1\endcsname . \quad}

\section{Fundamental group calculations (by Daniel Matei)}

If $K$ is a knot, we denote by $M_2(K)$ the $2$-fold cover of the 
$3$-sphere branched along $K$. It is well-known (e.g.~\cite{Ruberman,Viro}),
that the homeomorphism type of the closed $3$-manifold $M_2(K)$ is a
mutation invariant. Here we are mainly interested in its fundamental
group $\pi(K):=\pi_1(M_2(K))$. Let us review how a presentation of
$\pi(K)$ can be obtained from the knot group $G(K)$. 

The fundamental group $H(K)$ of the $2$-fold unbranched cover of the 
knot complement is the subgroup of $G(K)$ defined by the kernel of 
the homomorphism $G(K) \to \mathbb{Z}_2$ sending a meridian $\mu\in 
G(K)$ of $K$ (and in fact all meridians) to the generator of 
$\mathbb{Z}_2$. Then $\pi(K)$ is the quotient $H(K)/(\mu^2)$
of $H(K)$ by the subgroup normally generated by the squares of the
meridians of $K$, see~\cite{Rolfsen} for example. For our purposes
it is more convenient to view $\pi(K)$ as the index $2$ subgroup 
of $G(K)/(\mu^2)$ determined by the descending homomorphism 
$G(K)/(\mu^2) \to \mathbb{Z}_2$.

Suppose $K$ is given as the closure of the $n$-strand braid $b$,
written as a product in the standard generators $\sigma_k, 1\le k <n$,
and their inverses. We view braids as automorphisms of the free group 
on generators $x_1,\dots, x_n$, so that $\sigma_k(x_i)$ is equal to $x_k
x_{k+1} x_k^{-1}$ if $i=k$, to $x_k$ if $i=k+1$, and to $x_i$ otherwise.

We use the program GAP\ \cite{GAP}, and input braid generators and
braids as automorphisms of the free group like:

{\small
\begin{alltt}
\textit{gap>}f:=FreeGroup(4);
\textit{gap>}s1:=GroupHomomorphismByImages(f,f,[f.1,f.2,f.3,f.4],[f.1*f.2*f.1^-1,f.1,f.3,f.4]);
...
\textit{gap>}b1:=s1*s2*s2*s3^-1...
\end{alltt}
}

The knot group $G(K)$ has then a presentation with meridian generators $x_i$ 
and relations $b(x_i)=x_i$, where $1\le i\le n$. 
A presentation for $G(K)/(\mu^2)$  is obtained simply by adding $\mu^2=1$ to 
the above relations, with $\mu=x_i$ for some $i$. 

{\small
\begin{alltt}
\textit{gap>}rels:=List([1..4],i->Image(b1,GeneratorsOfGroup(f)[i])*GeneratorsOfGroup(f)[i]^-1);
\textit{gap>}rels:=Concatenation(rels,[GeneratorsOfGroup(f)[1]^2]);
\textit{gap>}q1:=f/rels;
\end{alltt}
}

{}From this presentation of $G(K)/(\mu^2)$, GAP then computes a
presentation for $\pi(K)$. First $\mathbb{Z}_2$ is realized as a
symmetric group of size $2$ and the homomorphism $G(K)/(\mu^2) \to
\mathbb{Z}_2$ is defined on generators by $x_i \to (1,2)$ for all
$i$, using the command \web|GroupHomomorphismByImages|. Then its
kernel $\pi(K)$ is found and converted into a finitely presented group. 

{\small
\begin{alltt}
\textit{gap>}c2:=Group((1,2));
\textit{gap>}p1:=Kernel(GroupHomomorphismByImages(q1,c2,GeneratorsOfGroup(q1),[(1,2),(1,2),(1,2),(1,2)])));
\textit{gap>}p1:=Image(IsomorphismFpGroup(p1));
\end{alltt}
}

Finally a presentation is created and displayed for the fundamental group 
of the branched cover.
For example, from the braid word in table \ref{tabg} for $14_{41763}$ we obtain (the
computer output is indicated by italics):

{\small
\begin{alltt}
\textit{gap>}TzPrint(PresentationFpGroup(p1))); \itshape
#I  generators: [ F1, F2, F3 ]
#I  relators:
#I  1.  11  [ 2, -3, 1, 3, -2, -2, 3, 1, -3, 2, -1 ]
#I  2.  12  [ 1, -2, 1, -2, 3, -1, 3, -2, 1, -2, 1, 2 ]
#I  3.  20  [ -3, 2, 2, -3, -3, 2, 2, -3, 2, 2, -3, -3, 2, 2, -3, 2, 1, -2, 1, 2 ]
\end{alltt}
}

For a pair of knots $K_1, K_2$ we will distinguish their double branched cover 
groups $\pi(K_1)$ and $\pi(K_2)$ using two types of numerical invariants.
Denote by $\pi$ either one of the two groups. The first invariant is simply 
the number of epimorphisms $\delta_{\Gamma}(\pi)$ of $\pi$ onto a finite group 
$\Gamma$ up to automorphisms of the target. The list of such epimorphisms is 
determined via the GAP command \verb|GQuotients(pi,Gamma)|. The second invariant 
is a list of abelianizations of certain finite index subgroups of $\pi$. The 
list is of two types: either a list $Ab_r(\pi)$ of the abelianizations of all 
(conjugacy classes of) index $r$ subgroups of $\pi$, or a list $Ab_\Gamma(\pi)$ 
of the abelianizations of the (conjugacy classes of) kernels of all epimorphisms 
from $\pi$ onto $\Gamma$. The finite index subgroups are determined via the GAP 
command \verb|LowIndexSubgroupsFpGroup(pi,TrivialSubgroup(pi),index)|. The kernels are obtained 
applying the command \verb|Kernel| to the homomorphisms in
\verb|GQuotients(pi,Gamma)|.
Finally, the abelianizations are obtained using the command 
\verb|AbelianInvariants|, and they are presented as lists of integers. For 
example, \verb|[ 0, 0, 2, 3, 3, 4 ]| stands for $\mathbb{Z}\times\mathbb{Z}
\times\mathbb{Z}_2\times\mathbb{Z}_3\times\mathbb{Z}_3\times\mathbb{Z}_4$.

As an example, let us verify that the homology group of $H(14_{41763})$
is what its order $\Dl_{14_{41763}}(-1)$, determined by the
Alexander polynomial $\Dl$, implies that it should be:

\def\trm#1{{\fontshape{n}\selectfont #1}}

{\small
\begin{alltt}
\textit{gap>}AbelianInvariants(p1);
\textit{[ 3 ]}
\end{alltt}
}

We start with three pairs $(K_1, K_2)$ of knots, whose groups $\pi(K_1)$ and 
$\pi(K_2)$ are not isomorphic, as their $\delta_{\Gamma}(\pi)$ are different: 
\begin{itemize}
\item for $(14_{41763},14_{42021})$ we have $\delta_{\Gamma}(\pi_1)=2$ 
and $\delta_{\Gamma}(\pi_2)=0$ with $\Gamma=Alt(7)$.

{\small
\begin{alltt}
\textit{gap>}Length(GQuotients(pi1,AlternatingGroup(7)));
\textit{2}
\textit{gap>}Length(GQuotients(pi2,AlternatingGroup(7)));
\textit{0}
\end{alltt}
}

\item for $(15_{219244},15_{228905})$ we have 
$\delta_{\Gamma}(\pi_1)=1$ and $\delta_{\Gamma}(\pi_2)=0$ with $\Gamma=PSL(2,13)$.

\item for $(15_{220504},15_{234873})$ we have
$\delta_{\Gamma}(\pi_1)=2$ and $\delta_{\Gamma}(\pi_2)=1$ with $\Gamma=PSL(2,7)$.

\end{itemize}

For the next three pairs we could not find (small enough) finite groups 
$\Gamma$ for which GAP returned (in a reasonable amount of time) 
different values of $\delta_{\Gamma}(\pi)$. We turned then to the 
abelianizations of finite index subgroups mentioned above.

For the pair $(14_{41739},14_{42126})$ it was possible to distinguish
$\pi(K_1)$ from $\pi(K_2)$ by considering the abelianizations
of the index $6$ subgroups.

{\small
\begin{alltt}\itshape
gap>\trm{lis:=LowIndexSubgroupsFpGroup(pi1,TrivialSubgroup( pi1 ),6);}
gap>\trm{lis6:=Filtered(lis,H->Index(pi1,H)=6);}
\end{alltt}
}

The second command is needed to select the subgroups of index $=6$; the
first command \webb o|LowIndexSubgroupsFpGroup| gives subgroups up to
given index. This feature depends on the algorithm GAP uses to find
the quotients, and is feasible only for subgroups of very small index.

Each group has $3$ index $6$ subgroups,
but the abelianizations of one of them on either side differ: 
% $Ab_6(\pi_1)=\{[9,9],[3,5,9],[2,2,16]\}$, 
% respectively $Ab_6(\pi_2)=\{[0,4,9],[3,5,9],[2,2,16]\}$.

{\small
\begin{alltt}\itshape
gap>\trm{Length(lis6)};
3
gap>\trm{List(lis6,H->AbelianInvariantsSubgroupFpGroup( pi1, H ));}
[ [ 9, 9 ], [ 3, 5, 9 ], [ 2, 2, 16 ] ]
...
gap>\trm{List(lis6,H->AbelianInvariantsSubgroupFpGroup( pi2, H ));}
[ [ 0, 4, 9 ], [ 3, 5, 9 ], [ 2, 2, 16 ] ]
\end{alltt}
}

For the next two pairs GAP could not determine enough low index 
subgroups of $\pi(K_1)$ and $\pi(K_2)$ in a reasonable amount of time. 
We then considered normal subgroups of higher index and computed their
abelianizations.

For the pair $(14_{41721},14_{42125})$ we first determined that
$\delta_{\Gamma}(\pi_1)=\delta_{\Gamma}(\pi_2)=1$ with $\Gamma=PSL(2,7)$.

{\small%
\begin{alltt}\itshape
gap>\trm{p1sl27:=GQuotients(p1,PSL(2,7));}
[ [ F1, F2, F3 ] -> [ (1,4)(2,5)(3,8)(6,7), (1,3,6,7)(2,4,5,8), (1,2,3,5)(4,8,7,6) ] ]
gap>\trm{p2sl27:=GQuotients(p2,PSL(2,7));}
[ [ F1, F2, F3, F4 ] -> [ (2,7,3,8,5,6,4), (1,2,4,7,3,5,8), (2,6,8,7,4,5,3), (2,3,5,4,7,8,6) ] ]
\end{alltt}
}

The two epimorphisms $h_1:\pi_1\to \Gamma$, respectively $h_2:\pi_2\to \Gamma$
are given by the images of the generators of $\pi_i$, treating $\Gamma=PSL(2,7)$ 
as a subgroup of $S_8$ and specifying the cycle decomposition of the permutations.

Then we computed using \verb|AbelianInvariants(Kernel(p1sl27[1]));| and
\webb l|AbelianInvariants(Kernel(p2sl27[1]));| 
the abelianizations $A_1$ and $A_2$ of the kernels of $h_{1,2}$. 
Since $A_1$ has $3$-torsion, whereas $A_2$ does not, we may conclude that 
$\pi(K_1)$ and $\pi(K_2)$ are not isomorphic.

Similarly, for the pair $(15_{148731},15_{156433})$ we first determined that 
$\delta_{\Gamma}(\pi_1)=\delta_{\Gamma}(\pi_2)=1$ with $\Gamma=Alt(6)$.
Then we found that the abelianizations $A_1$ and $A_2$ of the kernels of 
the two epimorphisms onto $Alt(6)$ are distinct, as $A_1$ is entirely torsion, 
whereas $A_2$ has free part of rank $10$. Thus $\pi(K_1)$ and $\pi(K_2)$ 
are not isomorphic. This last pair was the most difficult to break. The 
subgroups have a rather large index, $360$, the order of $Alt(6)$, and their
abelianizations are fairly complicated. 

Table \ref{tabg} summarizes the results.

% \newpage 

{

\parskip3mm

\newbox\@tempboxc
\setbox\@tempboxc=\hbox{\capt The data for the 6 pairs: knots, braid
  representation, number of generators and relations in the found
  presentation of $\pi$, and distinction method. The generators are
  numbered $f_1,\protect\dots,f_n$; in the relations we write $i$ for
  $f_i$ and use $\dl_i=f_if_{i+1}^{-1}f_i$ and $\gm_i=f_if_{i-1}^{-1}
  f_i$ to shorten the words. The abelianizations are abbreviated by
  omitting the product signs between the cyclic groups and writing
  $\bZ_*^k$ for (the product of) $k$ copies of $\bZ_*$.}
\def\x{\rule[-2mm]{\z@}{\z@}}
\captionwidth=0.7\vsize\relax
\@captionmargin0.08\vsize\relax

\begin{table}[ptb]
% \captionwidth0.8\vsize\relax
\newpage
\vbox to \textheight{\vfil
\tabcolsep5.8pt
\rottab{%
\def\hh{\\[2mm]\hline[1.5mm]}%
\def\hg{\\[2mm]\hline \hline [2mm]}%
\hbox to \textheight{\tiny\hss
\begin{mytab}{|c||p{5.5cm}||c|p{10cm}||p{3.5cm}|}%
  { &  & \multicolumn{2}{|r||}{ } & }%
  \hline [1.5mm]%
\mbox{pair} & \mbox{braids} & \mbox{\# gen of $\pi$} &
  \mbox{relators} & \mbox{distinction} \\[2mm]%
\hline
\hline [2mm]%
% table2.tex: this is a human edited copy of table.tex
$15_{219244}$ &
   { $\sg_{1} $ $ \sg_{2}^{2} $ $ \sg_{3}^{-1} $ $ \sg_{1} $ $ \sg_{2}^{-1} $ $ \sg_{1}^{2} $ $ \sg_{2}^{2} $ $ \sg_{3}^{-1} $ $ \sg_{2} $ $ \sg_{1} $ $ \sg_{2}^{-1} $ $ \sg_{3}^{-1} $ $ $} \x&
 3 &
    [ $\gm_{3} $ $ {3} $ $ \dl_{2}^{-1} $ $ {1} $ $ \dl_{1} $ $ {1} $ $ {2}^{-1} $ $ $];\ 
   [ ${1}^{-1} $ $ \dl_{1}^{-1} $ $ {2} $ $ {1}^{-2} $ $ {2} $ $ {3} $ $ {1}^{-1} $ $ {3} $ $ \gm_{3} $ $ {3} $ $ {1}^{-1} $ $ {3} $ $ {2} $ $ $];\ 
   [ ${2}^{-1} $ $ {1}^{2} $ $ {2}^{-1} $ $ {1}^{-1} $ $ \dl_{2} $ $ {1}^{-1} $ $ {2}^{-1} $ $ {1}^{2} $ $ {2}^{-1} $ $ {1}^{-1} $ $ \dl_{2} $ $ {1}^{-1} $ $ {2}^{-1} $ $ {1}^{2} $ $ {2}^{-1} $ $ {3}^{-1} $ $ $]\x &
$\delta_{PSL(2,13)}=1$ \hh
$15_{228905}$ &
   { $\sg_{1} $ $ \sg_{2}^{-1} $ $ \sg_{3}^{-1} $ $ \sg_{1} $ $ \sg_{2}^{2} $ $ \sg_{3}^{-1} $ $ \sg_{1}^{2} $ $ \sg_{2} $ $ \sg_{1}^{-1} $ $ \sg_{2} $ $ \sg_{3}^{-1} $ $ \sg_{2} $ $ \sg_{1} $ $ $} \x&
 3 &
    [ ${3}^{2} $ $ {2} $ $ {3} $ $ {2} $ $ {1}^{-1} $ $ {2} $ $ {3} $ $ {2} $ $ {3}^{2} $ $ {1}^{-2} $ $ $];\ 
   [ ${2} $ $ {1}^{-3} $ $ {2} $ $ {3}^{-1} $ $ {1}^{3} $ $ {2}^{-1} $ $ {3}^{-1} $ $ {2}^{-1} $ $ {3}^{-1} $ $ {2}^{-1} $ $ {1}^{3} $ $ {3}^{-1} $ $ $];\ 
   [ ${1}^{3} $ $ {2}^{-1} $ $ {3} $ $ {1}^{-2} $ $ {3}^{2} $ $ {1}^{-2} $ $ {3} $ $ {2}^{-1} $ $ {1}^{3} $ $ \dl_{2}^{-1} $ $ $]\x &
$\delta_{PSL(2,13)}=0$\hg
$14_{41763}$ &
   { $\sg_{1}^{-1} $ $ \sg_{2}^{-1} $ $ \sg_{3}^{-1} $ $ \sg_{4} $ $ \sg_{3}^{-1} $ $ \sg_{4} $ $ \sg_{2}^{-1} $ $ \sg_{4} $ $ \sg_{1} $ $ \sg_{3}^{-1} $ $ \sg_{2} $ $ \sg_{3}^{-1} $ $ \sg_{4} $ $ \sg_{3}^{2} $ $ \sg_{2} $ $ $} \x&
 3 &
    [ ${2} $ $ {3}^{-1} $ $ {1} $ $ {3} $ $ {2}^{-2} $ $ {3} $ $ {1} $ $ {3}^{-1} $ $ {2} $ $ {1}^{-1} $ $ $];\ 
   [ $\dl_{1} $ $ {2}^{-1} $ $ {3} $ $ {1}^{-1} $ $ {3} $ $ {2}^{-1} $ $ \dl_{1} $ $ {2} $ $ $];\ 
   [ ${3}^{-1} $ $ {2}^{2} $ $ {3}^{-2} $ $ {2} $ $ \dl_{2} $ $ {2} $ $ {3}^{-2} $ $ {2} $ $ \dl_{2} $ $ \dl_{1} $ $ {2} $ $ $]\x &
$\delta_{Alt(7)}=2$  \hh
$14_{42021}$ &
   { $\sg_{1}^{-1} $ $ \sg_{2} $ $ \sg_{1}^{-1} $ $ \sg_{2} $ $ \sg_{3} $ $ \sg_{4} $ $ \sg_{3}^{-1} $ $ \sg_{2}^{-2} $ $ \sg_{1}^{-1} $ $ \sg_{2} $ $ \sg_{1}^{-1} $ $ \sg_{3}^{-1} $ $ \sg_{2} $ $ \sg_{4}^{-1} $ $ \sg_{3} $ $ $} \x&
 3 &
    [ ${1}^{-1} $ $ \gm_{3}^{-1} $ $ {1}^{-1} $ $ \dl_{1}^{-1} $ $ {2} $ $ {1}^{-1} $ $ $];\ 
   [ ${1}^{-2} $ $ {2} $ $ {3}^{-1} $ $ \dl_{1} $ $ {1} $ $ \dl_{1} $ $ {3}^{-1} $ $ {2} $ $ {1}^{-2} $ $ {3} $ $ $];\ 
   [ ${1} $ $ {2}^{-1} $ $ {3} $ $ \dl_{1}^{-1} $ $ \gm_{3}^{2} $ $ \dl_{1}^{-1} $ $ \gm_{3}^{2} $ $ \dl_{1}^{-1} $ $ {3} $ $ {2}^{-1} $ $ {1} $ $ $]\x &
$\delta_{Alt(7)}=0$ \hg
$15_{148731}$ &
   { $\sg_{1}^{-1} $ $ \sg_{2}^{-1} $ $ \sg_{3}^{-1} $ $ \sg_{4}^{2} $ $ \sg_{5}^{-1} $ $ \sg_{4} $ $ \sg_{3}^{-1} $ $ \sg_{2} $ $ \sg_{1} $ $ \sg_{3} $ $ \sg_{4} $ $ \sg_{3} $ $ \sg_{2} $ $ \sg_{3} $ $ \sg_{4}^{-2} $ $ \sg_{5} $ $ \sg_{4}^{-1} $ $ \sg_{3} $ $ \sg_{2}^{-1} $ $ \sg_{4}^{-1} $ $ \sg_{3} $ $ $} \x&
 3 &
    [ ${2}^{2} $ $ \dl_{1} $ $ {2}^{2} $ $ {3}^{2} $ $ $];\ 
   [ ${1}^{-1} $ $ {2} $ $ {3} $ $ {1} $ $ {3}^{-1} $ $ {2}^{-2} $ $ {3}^{-1} $ $ {1} $ $ {3} $ $ {2} $ $ $];\ 
   [ ${2}^{-1} $ $ {3}^{-1} $ $ {1}^{-1} $ $ {3}^{-1} $ $ {2}^{-1} $ $ \dl_{1} $ $ {3}^{-1} $ $ \dl_{1} $ $ {2}^{-1} $ $ {3}^{-1} $ $ {1}^{-1} $ $ {3}^{-1} $ $ {2}^{-2} $ $ $]\x &
$Ab_{Alt(6)}=$ [$ \x\bZ_{2}^{41} \bZ_{3}^{17} \bZ_{4}^{60}$ $ \bZ_{5}^{10} \bZ_{9}^{7} \bZ_{19}^{5}$] \hh
$15_{156433}$ &
   { $\sg_{1}^{-1} $ $ \sg_{2}^{-1} $ $ \sg_{3}^{-1} $ $ \sg_{4} $ $ \sg_{5} $ $ \sg_{3}^{-1} $ $ \sg_{2}^{2} $ $ \sg_{3}^{-1} $ $ \sg_{4}^{-1} $ $ \sg_{3} $ $ \sg_{2}^{-1} $ $ \sg_{1} $ $ \sg_{2}^{-2} $ $ \sg_{3} $ $ \sg_{4} $ $ \sg_{3}^{-1} $ $ \sg_{5}^{-1} $ $ \sg_{4}^{-1} $ $ \sg_{2} $ $ \sg_{3}^{-1} $ $ \sg_{2} $ $ $} \x&
 4 &
    [ $\dl_{1}^{-1} $ $ \dl_{3} $ $ {4} $ $ {3}^{-1} $ $ \dl_{1} $ $ {3}^{-1} $ $ \dl_{1} $ $ {3}^{-1} $ $ {4} $ $ \dl_{3} $ $ {1}^{-1} $ $ {2} $ $ $];\ 
   [ ${4} $ $ {3}^{-1} $ $ \dl_{1} $ $ {3}^{-1} $ $ \dl_{1} $ $ \dl_{3}^{-1} $ $ {4} $ $ {3}^{-1} $ $ \dl_{1} $ $ {2}^{-1} $ $ {1} $ $ \dl_{3}^{-1} $ $ {1} $ $ $];\ 
   [ $\dl_{1} $ $ {3}^{-1} $ $ {2} $ $ \dl_{3}^{-1} $ $ \dl_{1} $ $ {2}^{-1} $ $ {1} $ $ \dl_{3}^{-2} $ $ \dl_{1} $ $ {2} $ $ {1}^{-1} $ $ \dl_{3} $ $ {1}^{-1} $ $ {2} $ $ {3}^{-1} $ $ \dl_{1} $ $ \dl_{3}^{-1} $ $ $];\ 
   [ ${2} $ $ \dl_{3}^{-1} $ $ \dl_{1} $ $ {2}^{-1} $ $ {1} $ $ \dl_{3}^{-2} $ $ \dl_{1} $ $ {4}^{-1} $ $ {3} $ $ \dl_{1}^{-1} $ $ {3} $ $ {2}^{-1} $ $ {1} $ $ \dl_{3}^{-1} $ $ {1} $ $ {2}^{-1} $ $ {4}^{-1} $ $ {3} $ $ \dl_{1}^{-1} $ $ {3} $ $ {2}^{-1} $ $ {1} $ $ \dl_{3}^{-1} $ $ {1} $ $ $]\x &
$Ab_{Alt(6)}=$ [$ \x\bZ^{10} \bZ_{2}^{17} \bZ_{3}^{21} \bZ_{5}^{10}$ $ \bZ_{8}^{16} \bZ_{16} \bZ_{17}^{5} \bZ_{32}^{5} \bZ_{64}^{4}$]\hg
$15_{220504}$ &
   { $\sg_{1}^{-1} $ $ \sg_{2} $ $ \sg_{3}^{-1} $ $ \sg_{4} $ $ \sg_{5} $ $ \sg_{2} $ $ \sg_{3}^{2} $ $ \sg_{2} $ $ \sg_{1} $ $ \sg_{2}^{2} $ $ \sg_{3}^{-1} $ $ \sg_{4}^{-1} $ $ \sg_{2} $ $ \sg_{3}^{-1} $ $ \sg_{2} $ $ \sg_{3}^{-1} $ $ \sg_{4}^{-1} $ $ \sg_{5}^{-1} $ $ \sg_{4}^{-1} $ $ \sg_{3} $ $ \sg_{2} $ $ $} \x&
 4 &
    [ ${3} $ $ {1}^{-1} $ $ \dl_{1}^{-1} $ $ {1}^{-1} $ $ \dl_{1}^{-1} $ $ {1}^{-1} $ $ {3} $ $ $];\ 
   [ ${4}^{-1} $ $ {1}^{3} $ $ {4}^{-1} $ $ {1}^{3} $ $ {4}^{-1} $ $ {3} $ $ {4}^{-1} $ $ {3}^{-1} $ $ {4}^{-1} $ $ {3} $ $ $];\ 
   [ ${4}^{-1} $ $ {3}^{-1} $ $ {1} $ $ \dl_{1} $ $ {4} $ $ {1}^{-3} $ $ {4} $ $ \dl_{1} $ $ {1} $ $ {3}^{-1} $ $ $];\ 
   [ ${3} $ $ {4} $ $ {3}^{-1} $ $ {4} $ $ {3}^{2} $ $ {1}^{-2} $ $ {2}^{2} $ $ {1}^{-2} $ $ {3}^{2} $ $ {1} $ $ {2} $ $ {1}^{-2} $ $ {3} $ $ $]\x &
 $\delta_{PSL(2,7)}=2$ \hh
$15_{234873}$ &
   { $\sg_{1}^{-1} $ $ \sg_{2} $ $ \sg_{3}^{-1} $ $ \sg_{4}^{-1} $ $ \sg_{5} $ $ \sg_{2} $ $ \sg_{3}^{-1} $ $ \sg_{4}^{-1} $ $ \sg_{2}^{2} $ $ \sg_{3}^{-1} $ $ \sg_{2} $ $ \sg_{1} $ $ \sg_{2} $ $ \sg_{3}^{2} $ $ \sg_{4}^{-1} $ $ \sg_{5}^{-1} $ $ \sg_{2} $ $ \sg_{3}^{-1} $ $ \sg_{4} $ $ \sg_{3} $ $ \sg_{2} $ $ $} \x&
 3 &
    [ ${2}^{-1} $ $ {3}^{-1} $ $ {1}^{-2} $ $ {3}^{-1} $ $ {2}^{-2} $ $ {3}^{-1} $ $ {1}^{-2} $ $ {3}^{-1} $ $ {2}^{-1} $ $ {3} $ $ {1}^{-1} $ $ {3} $ $ $];\ 
   [ $\dl_{2} $ $ {3} $ $ {1}^{2} $ $ {3} $ $ {2} $ $ {3} $ $ {1}^{2} $ $ {3} $ $ \dl_{2} $ $ {1}^{-3} $ $ $];\ 
   [ ${1}^{3} $ $ {3} $ $ {1}^{2} $ $ {3} $ $ {2}^{-1} $ $ {3}^{-1} $ $ {1}^{-2} $ $ {3}^{-1} $ $ {2}^{-1} $ $ {3} $ $ {1}^{2} $ $ \dl_{2}^{-1} $ $ {1}^{2} $ $ \dl_{1} $ $ {1}^{2} $ $ {3} $ $ {1}^{3} $ $ {2}^{-1} $ $ $]\x &
$\delta_{PSL(2,7)}=1$ \hg
$14_{41721}$ &
   { $\sg_{1}^{-1} $ $ \sg_{2} $ $ \sg_{3} $ $ \sg_{4}^{-1} $ $ \sg_{5}^{-1} $ $ \sg_{6}^{-1} $ $ \sg_{5} $ $ \sg_{4}^{-1} $ $ \sg_{3}^{-1} $ $ \sg_{2}^{-1} $ $ \sg_{4} $ $ \sg_{1} $ $ \sg_{3} $ $ \sg_{2} $ $ \sg_{4}^{-1} $ $ \sg_{3} $ $ \sg_{4} $ $ \sg_{5}^{-1} $ $ \sg_{4} $ $ \sg_{3}^{-1} $ $ \sg_{2}^{-1} $ $ \sg_{6} $ $ \sg_{5}^{-1} $ $ \sg_{4}^{-1} $ $ \sg_{3}^{-1} $ $ \sg_{4} $ $ \sg_{5} $ $ \sg_{4}^{-1} $ $ \sg_{5} $ $ \sg_{4} $ $ $} \x&
 3 &
    [ ${3} $ $ {1}^{-2} $ $ \dl_{2} $ $ {1}^{-2} $ $ {2} $ $ \dl_{1} $ $ {2}^{-1} $ $ {1} $ $ {2} $ $ {1}^{-2} $ $ \dl_{2} $ $ {1}^{-2} $ $ $];\ 
   [ ${3}^{-1} $ $ {2} $ $ {1}^{-1} $ $ \dl_{1}^{-1} $ $ {2} $ $ {1}^{-1} $ $ \dl_{1}^{-1} $ $ {2} $ $ {1}^{-2} $ $ \dl_{1}^{-1} $ $ {2} $ $ {1}^{-1} $ $ \dl_{1}^{-1} $ $ {2} $ $ {1}^{-2} $ $ \dl_{2} $ $ {1}^{-2} $ $ {2} $ $ $];\ 
   [ ${3} $ $ {2}^{-1} $ $ {1} $ $ \dl_{1} $ $ {2}^{-1} $ $ {1} $ $ \dl_{1} $ $ {2}^{-1} $ $ {1}^{2} $ $ {2}^{-1} $ $ {3} $ $ {1}^{-2} $ $ \dl_{2} $ $ {1}^{-1} $ $ \dl_{1}^{-1} $ $ {2} $ $ \dl_{1}^{-1} $ $ {1}^{-1} $ $ \dl_{2} $ $ {1}^{-2} $ $ $];\ 
   [ $\dl_{1}^{2} $ $ {2}^{-1} $ $ {1}^{2} $ $ \dl_{2}^{-1} $ $ {3} $ $ {2}^{-1} $ $ {1} $ $ \dl_{1} $ $ {2}^{-1} $ $ {1} $ $ \dl_{1} $ $ {2}^{-1} $ $ {1} $ $ \dl_{1} $ $ {1} $ $ \dl_{2}^{-1} $ $ {1} $ $ \dl_{1} $ $ {1} $ $ {2}^{-1} $ $ $]\x &
 $Ab_{PSL(2,7)}=$ [ $\x\bZ_{2} $ $ \bZ_{3}^{15} $ $ \bZ_{5}^{19} $ $ \bZ_{7}^{3} $ $ \bZ_{11}^{8} $ $ $ ] \hh
$14_{42125}$ &
   { $\sg_{1}^{-1} $ $ \sg_{2} $ $ \sg_{3} $ $ \sg_{4} $ $ \sg_{5}^{-1} $ $ \sg_{6}^{-1} $ $ \sg_{5} $ $ \sg_{4} $ $ \sg_{3}^{-1} $ $ \sg_{2}^{-1} $ $ \sg_{4}^{-1} $ $ \sg_{1} $ $ \sg_{3} $ $ \sg_{2} $ $ \sg_{4} $ $ \sg_{3} $ $ \sg_{4}^{-1} $ $ \sg_{5}^{-1} $ $ \sg_{4}^{-1} $ $ \sg_{3}^{-1} $ $ \sg_{2}^{-1} $ $ \sg_{6} $ $ \sg_{4} $ $ \sg_{5}^{-1} $ $ \sg_{4} $ $ \sg_{3}^{-1} $ $ \sg_{4}^{-1} $ $ \sg_{5}^{2} $ $ \sg_{4}^{-1} $ $ $} \x&
 4 &
    [ ${1}^{-1} $ $ {2} $ $ {4}^{-1} $ $ {2} $ $ {4}^{-1} $ $ {2} $ $ {1}^{-1} $ $ $];\ 
   [ $\dl_{3}^{-1} $ $ {2}^{2} $ $ {1}^{-1} $ $ {4} $ $ {2}^{-1} $ $ {1} $ $ \dl_{2}^{-1} $ $ {1} $ $ {2}^{-1} $ $ {4} $ $ $];\ 
   [ ${1} $ $ \dl_{2} $ $ {1} $ $ {2}^{-1} $ $ {4}^{2} $ $ {2}^{-1} $ $ {1} $ $ \dl_{2} $ $ {1}^{-1} $ $ {2} $ $ {4}^{-1} $ $ {1}^{-1} $ $ {2} $ $ {4}^{-1} $ $ $];\ 
   [ ${4} $ $ {2}^{-1} $ $ {1} $ $ \dl_{2}^{-1} $ $ {1}^{-1} $ $ {2} $ $ {4}^{-1} $ $ {2}^{-1} $ $ {4} $ $ {3}^{-1} $ $ {2} $ $ {4} $ $ {2}^{-1} $ $ {1} $ $ \dl_{2}^{-1} $ $ {1}^{-1} $ $ {2} $ $ {4}^{-1} $ $ {2}^{-1} $ $ {3} $ $ {4}^{-1} $ $ {2} $ $ {1}^{-1} $ $ \dl_{2} $ $ $]\x &
$Ab_{PSL(2,7)}=$ [ $\x\bZ_{2} $ $ \bZ_{7}^{3} $ $ \bZ_{17}^{7} $ $ \bZ_{89}^{6} $ $ $] \hg
$14_{41739}$ &
   { $\sg_{1}^{-1} $ $ \sg_{2} $ $ \sg_{3} $ $ \sg_{4}^{-1} $ $ \sg_{5}^{-1} $ $ \sg_{6}^{-1} $ $ \sg_{5} $ $ \sg_{4}^{-1} $ $ \sg_{3}^{-1} $ $ \sg_{2}^{-1} $ $ \sg_{4} $ $ \sg_{1} $ $ \sg_{3} $ $ \sg_{2} $ $ \sg_{4}^{-1} $ $ \sg_{3} $ $ \sg_{4} $ $ \sg_{5}^{-1} $ $ \sg_{4} $ $ \sg_{3}^{-1} $ $ \sg_{2}^{-1} $ $ \sg_{6} $ $ \sg_{5}^{-1} $ $ \sg_{4}^{-1} $ $ \sg_{3}^{-1} $ $ \sg_{4} $ $ \sg_{5}^{-2} $ $ \sg_{4} $ $ \sg_{5}^{-1} $ $ $} \x&
 3 &
    [ ${3}^{-1} $ $ {2} $ $ \dl_{1} $ $ {2}^{-1} $ $ {1} $ $ \dl_{2} $ $ {3}^{-1} $ $ {2} $ $ {1}^{-1} $ $ \dl_{2} $ $ {3}^{-1} $ $ {2} $ $ \dl_{1} $ $ {2}^{-1} $ $ {1} $ $ {2} $ $ {3}^{-1} $ $ {1}^{2} $ $ $];\ 
   [ $\dl_{1}^{-1} $ $ {2} $ $ {1}^{-1} $ $ \dl_{1}^{-2} $ $ {2} $ $ {1}^{-1} $ $ \dl_{2}^{-1} $ $ {3} $ $ {2}^{-1} $ $ \dl_{1}^{-1} $ $ {2} $ $ {1}^{-1} $ $ \dl_{2}^{-1} $ $ \dl_{1}^{-1} $ $ {2} $ $ {1}^{-1} $ $ \dl_{2}^{-1} $ $ {3} $ $ {2}^{-1} $ $ $];\ 
   [ ${2} $ $ {1}^{-2} $ $ \dl_{2} $ $ {3}^{-1} $ $ {2} $ $ \dl_{1} $ $ {2}^{-1} $ $ {1} $ $ \dl_{2} $ $ {3}^{-1} $ $ {2} $ $ {1}^{-2} $ $ {2} $ $ \dl_{1} $ $ {2}^{-1} $ $ {1} $ $ \dl_{2} $ $ {3}^{-1} $ $ {2} $ $ \dl_{1} $ $ {2}^{-1} $ $ {1} $ $ $]\x &
 $Ab_{6}=$ [ $[ \x\bZ_{9}^{2} $ $ $], [ $\bZ_{3} $ $ \bZ_{5} $ $ \bZ_{9} $ $ $], [ $\bZ_{2}^{2} $ $ \bZ_{16} $ $ $] ]  \hh
$14_{42126}$ &
   { $\sg_{1}^{-1} $ $ \sg_{2} $ $ \sg_{3} $ $ \sg_{4} $ $ \sg_{5}^{-1} $ $ \sg_{6}^{-1} $ $ \sg_{5} $ $ \sg_{4} $ $ \sg_{3}^{-1} $ $ \sg_{2}^{-1} $ $ \sg_{4}^{-1} $ $ \sg_{1} $ $ \sg_{3} $ $ \sg_{2} $ $ \sg_{4} $ $ \sg_{3} $ $ \sg_{4}^{-1} $ $ \sg_{5}^{-1} $ $ \sg_{4}^{-1} $ $ \sg_{3}^{-1} $ $ \sg_{2}^{-1} $ $ \sg_{6} $ $ \sg_{4} $ $ \sg_{5}^{-1} $ $ \sg_{4} $ $ \sg_{3}^{-1} $ $ \sg_{4}^{-1} $ $ \sg_{5}^{-2} $ $ \sg_{4}^{-1} $ $ $} \x&
 3 &
    [ ${1} $ $ \dl_{2}^{-1} $ $ {1}^{2} $ $ \dl_{2}^{-1} $ $ {1}^{4} $ $ \dl_{2}^{-1} $ $ {1}^{2} $ $ \dl_{2}^{-1} $ $ {1} $ $ {3} $ $ {1}^{-1} $ $ {3} $ $ $];\ 
   [ ${3} $ $ {2}^{-1} $ $ {1} $ $ {3} $ $ {1}^{-1} $ $ {3} $ $ {2}^{-1} $ $ {1}^{-1} $ $ \dl_{2} $ $ {1}^{-2} $ $ \dl_{2} $ $ {1}^{-1} $ $ {2}^{-1} $ $ {3} $ $ {1}^{-1} $ $ {3} $ $ {1} $ $ \dl_{2}^{-1} $ $ {1} $ $ \dl_{1} $ $ {1} $ $ {2}^{-1} $ $ $];\ 
   [ ${1}^{-2} $ $ {3} $ $ {2}^{-1} $ $ {1}^{-1} $ $ \dl_{2} $ $ {1}^{-2} $ $ \dl_{2} $ $ {1}^{-1} $ $ {2}^{-1} $ $ {1}^{-1} $ $ \dl_{2} $ $ {1}^{-2} $ $ \dl_{2} $ $ {1}^{-1} $ $ {2}^{-1} $ $ {3} $ $ {1}^{-2} $ $ \dl_{2} $ $ $]\x &
 $Ab_{6}=$ [ $[ \x\bZ $ $ \bZ_{4} $ $ \bZ_{9} $ $ $], [ $\bZ_{3} $ $ \bZ_{5} $ $ \bZ_{9} $ $ $], [ $\bZ_{2}^{2} $ $ \bZ_{16} $ $ $] ] \hg

% \\[2mm]
% \hline
% \hline [2mm]%
% total & & & & \\[1.2mm]
% \hline
\end{mytab}%
\hss}%
}{\unhcopy\@tempboxc
  \protect\label{tabg}}
\vss}
\newpage
\setbox\@tempboxc=\null
\end{table}
}

{\small\it Daniel Matei, 
Institute of Mathematics "Simion Stoilow" of the Romanian Academy
P.O. BOX 1-764, RO-014700 Bucharest, Romania; daniel.matei@imar.ro  \\
and \\
Graduate School of Mathematical Sciences, University of Tokyo,
3-8-1 Komaba Meguro-ku, Tokyo 153-8914, Japan
} \\
% }
\end{appendix}

\end{document}